\documentclass[11pt]{amsart}
\usepackage{amsmath}
\usepackage{amssymb}
\usepackage{amsfonts}
\usepackage{amsthm}
\usepackage{enumerate}
\usepackage[mathscr]{eucal}
\usepackage{verbatim}

\begin{document}

\newtheorem{theorem}{Theorem} [section]
\newtheorem{proposition}[theorem]{Proposition}
\newtheorem{conjecture}[theorem]{Conjecture}
\def\theconjecture{\unskip}
\newtheorem{corollary}[theorem]{Corollary}
\newtheorem{lemma}[theorem]{Lemma}
\newtheorem{observation}[theorem]{Observation}
\theoremstyle{definition}
\newtheorem{definition}{Definition}
\newtheorem{remark}{Remark}
\def\theremark{\unskip}
\newtheorem{question}{Question}
\def\thequestion{\unskip}
\newtheorem{example}{Example}
\def\theexample{\unskip}
\newtheorem{problem}{Problem}

\numberwithin{theorem}{section}
\numberwithin{definition}{section}
\numberwithin{equation}{section}

\theoremstyle{plain}
\newtheorem{thmsub}{Theorem}[subsection]
\newtheorem{lemmasub}[thmsub]{Lemma}
\newtheorem{corollarysub}[thmsub]{Corollary}
\newtheorem{propositionsub}[thmsub]{Proposition}
\newtheorem{defnsub}[thmsub]{Definition}
\numberwithin{equation}{section}

\def\prob{\mu}
\def\fpqs{{F^{pq}_{\sigma}}}
\def\intslash{\rlap{\kern  .32em $\mspace {.5mu}\backslash$ }\int}
\def\qsl{{\rlap{\kern  .32em $\mspace {.5mu}\backslash$ }\int_{Q}}}
\def\Re{\operatorname{Re\,}}
\def\Im{\operatorname{Im\,}}
\def\mx{{\max}}
\def\mn{{\min}}
\def\vth{\vartheta}
\def\tDel{{\widetilde \Delta}}
\def\eann{{{\mathcal E_o}}}
\def\rn{\rr^{n}}
\def\rr{\mathbb R}
\def\R{\mathbb R}
\def\Q{\mathcal Q}
\def\N{\mathbb N}
\def\complex{{\mathbb C}}
\def\norm#1{{ \left|  #1 \right| }}
\def\Norm#1{{ \left\|  #1 \right\| }}
\def\set#1{{ \left\{ #1 \right\} }}
\def\floor#1{{\lfloor #1 \rfloor }}
\def\emph#1{{\it #1 }}
\def\diam{{\text{\rm diam}}}
\def\osc{{\text{\rm osc}}}
\def\ffB{\mathcal B}
\def\itemize#1{\item"{#1}"}
\def\seq{\subseteq}
\def\Id{\text{\sl Id}}

\def\Ga{\Gamma}
\def\ga{\gamma}
\def\Th{\Theta}

\def\prd{{\text{\it prod}}}
\def\parab{{\text{\it parabolic}}}

\def\eg{{\it e.g. }}
\def\cf{{\it cf}}
\def\Rn{{\mathbb R^n}}
\def\Rd{{\mathbb R^d}}
\def\sgn{{\text{\rm sign }}}
\def\rank{{\text{\rm rank }}}
\def\corank{{\text{\rm corank }}}
\def\coker{{\text{\rm Coker }}}
\def\loc{{\text{\rm loc}}}
\def\spec{{\text{\rm spec}}}

\def\comp{{\text{\rm comp}}}

\def\Coi{{C^\infty_0}}
\def\dist{{\text{\rm dist}}}
\def\diag{{\text{\rm diag}}}
\def\supp{{\text{\rm supp }}}
\def\rad{{\text{\rm rad}}}
\def\Lip{{\text{\rm Lip}}}
\def\inn#1#2{\langle#1,#2\rangle}
\def\biginn#1#2{\big\langle#1,#2\big\rangle}
\def\rta{\rightarrow}
\def\lta{\leftarrow}
\def\noi{\noindent}
\def\lcontr{\rfloor}
\def\lco#1#2{{#1}\lcontr{#2}}
\def\lcoi#1#2{\imath({#1}){#2}}
\def\rco#1#2{{#1}\rcontr{#2}}
\def\bin#1#2{{\pmatrix {#1}\\{#2}\endpmatrix}}
\def\meas{{\text{\rm meas}}}

\def\card{\text{\rm card}}
\def\lc{\lesssim}
\def\gc{\gtrsim}
\def\pv{\text{\rm p.v.}}

\def\alp{\alpha}             \def\Alp{\Alpha}
\def\bet{\beta}
\def\gam{\gamma}             \def\Gam{\Gamma}
\def\del{\delta}             \def\Del{\Delta}
\def\eps{\varepsilon}
\def\ep{\epsilon}
\def\zet{\zeta}
\def\tet{\theta}             \def\Tet{\Theta}
\def\iot{\iota}
\def\kap{\kappa}
\def\ka{\kappa}
\def\lam{\lambda}            \def\Lam{\Lambda}
\def\la{\lambda}             \def\La{\Lambda}
\def\sig{\sigma}             \def\Sig{\Sigma}
\def\si{\sigma}              \def\Si{\Sigma}
\def\vphi{\varphi}
\def\ome{\omega}             \def\Ome{\Omega}
\def\om{\omega}              \def\Om{\Omega}

\def\fA{{\mathfrak {A}}}
\def\fB{{\mathfrak {B}}}
\def\fC{{\mathfrak {C}}}
\def\fD{{\mathfrak {D}}}
\def\fE{{\mathfrak {E}}}
\def\fF{{\mathfrak {F}}}
\def\fG{{\mathfrak {G}}}
\def\fH{{\mathfrak {H}}}
\def\fI{{\mathfrak {I}}}
\def\fJ{{\mathfrak {J}}}
\def\fK{{\mathfrak {K}}}
\def\fL{{\mathfrak {L}}}
\def\fM{{\mathfrak {M}}}
\def\fN{{\mathfrak {N}}}
\def\fO{{\mathfrak {O}}}
\def\fP{{\mathfrak {P}}}
\def\fQ{{\mathfrak {Q}}}
\def\fR{{\mathfrak {R}}}
\def\fS{{\mathfrak {S}}}
\def\fT{{\mathfrak {T}}}
\def\fU{{\mathfrak {U}}}
\def\fV{{\mathfrak {V}}}
\def\fW{{\mathfrak {W}}}
\def\fX{{\mathfrak {X}}}
\def\fY{{\mathfrak {Y}}}
\def\fZ{{\mathfrak {Z}}}

\def\fa{{\mathfrak {a}}}
\def\fb{{\mathfrak {b}}}
\def\fc{{\mathfrak {c}}}
\def\fd{{\mathfrak {d}}}
\def\fe{{\mathfrak {e}}}
\def\ff{{\mathfrak {f}}}
\def\fg{{\mathfrak {g}}}
\def\fh{{\mathfrak {h}}}
\def\fj{{\mathfrak {j}}}
\def\fk{{\mathfrak {k}}}
\def\fl{{\mathfrak {l}}}
\def\fm{{\mathfrak {m}}}
\def\fn{{\mathfrak {n}}}
\def\fo{{\mathfrak {o}}}
\def\fp{{\mathfrak {p}}}
\def\fq{{\mathfrak {q}}}
\def\fr{{\mathfrak {r}}}
\def\fs{{\mathfrak {s}}}
\def\ft{{\mathfrak {t}}}
\def\fu{{\mathfrak {u}}}
\def\fv{{\mathfrak {v}}}
\def\fw{{\mathfrak {w}}}
\def\fx{{\mathfrak {x}}}
\def\fy{{\mathfrak {y}}}
\def\fz{{\mathfrak {z}}}

\def\bbA{{\mathbb {A}}}
\def\bbB{{\mathbb {B}}}
\def\bbC{{\mathbb {C}}}
\def\bbD{{\mathbb {D}}}
\def\bbE{{\mathbb {E}}}
\def\bbF{{\mathbb {F}}}
\def\bbG{{\mathbb {G}}}
\def\bbH{{\mathbb {H}}}
\def\bbI{{\mathbb {I}}}
\def\bbJ{{\mathbb {J}}}
\def\bbK{{\mathbb {K}}}
\def\bbL{{\mathbb {L}}}
\def\bbM{{\mathbb {M}}}
\def\bbN{{\mathbb {N}}}
\def\bbO{{\mathbb {O}}}
\def\bbP{{\mathbb {P}}}
\def\bbQ{{\mathbb {Q}}}
\def\bbR{{\mathbb {R}}}
\def\bbS{{\mathbb {S}}}
\def\bbT{{\mathbb {T}}}
\def\bbU{{\mathbb {U}}}
\def\bbV{{\mathbb {V}}}
\def\bbW{{\mathbb {W}}}
\def\bbX{{\mathbb {X}}}
\def\bbY{{\mathbb {Y}}}
\def\bbZ{{\mathbb {Z}}}

\def\cA{{\mathcal {A}}}
\def\cB{{\mathcal {B}}}
\def\cC{{\mathcal {C}}}
\def\cD{{\mathcal {D}}}
\def\cE{{\mathcal {E}}}
\def\cF{{\mathcal {F}}}
\def\cG{{\mathcal {G}}}
\def\cH{{\mathcal {H}}}
\def\cI{{\mathcal {I}}}
\def\cJ{{\mathcal {J}}}
\def\cK{{\mathcal {K}}}
\def\cL{{\mathcal {L}}}
\def\cM{{\mathcal {M}}}
\def\cN{{\mathcal {N}}}
\def\cO{{\mathcal {O}}}
\def\cP{{\mathcal {P}}}
\def\cQ{{\mathcal {Q}}}
\def\cR{{\mathcal {R}}}
\def\cS{{\mathcal {S}}}
\def\cT{{\mathcal {T}}}
\def\cU{{\mathcal {U}}}
\def\cV{{\mathcal {V}}}
\def\cW{{\mathcal {W}}}
\def\cX{{\mathcal {X}}}
\def\cY{{\mathcal {Y}}}
\def\cZ{{\mathcal {Z}}}

\def\tA{{\widetilde{A}}}
\def\tB{{\widetilde{B}}}
\def\tC{{\widetilde{C}}}
\def\tD{{\widetilde{D}}}
\def\tE{{\widetilde{E}}}
\def\tF{{\widetilde{F}}}
\def\tG{{\widetilde{G}}}
\def\tH{{\widetilde{H}}}
\def\tI{{\widetilde{I}}}
\def\tJ{{\widetilde{J}}}
\def\tK{{\widetilde{K}}}
\def\tL{{\widetilde{L}}}
\def\tM{{\widetilde{M}}}
\def\tN{{\widetilde{N}}}
\def\tO{{\widetilde{O}}}
\def\tP{{\widetilde{P}}}
\def\tQ{{\widetilde{Q}}}
\def\tR{{\widetilde{R}}}
\def\tS{{\widetilde{S}}}
\def\tT{{\widetilde{T}}}
\def\tU{{\widetilde{U}}}
\def\tV{{\widetilde{V}}}
\def\tW{{\widetilde{W}}}
\def\tX{{\widetilde{X}}}
\def\tY{{\widetilde{Y}}}
\def\tZ{{\widetilde{Z}}}

\def\tcA{{\widetilde{\mathcal {A}}}}
\def\tcB{{\widetilde{\mathcal {B}}}}
\def\tcC{{\widetilde{\mathcal {C}}}}
\def\tcD{{\widetilde{\mathcal {D}}}}
\def\tcE{{\widetilde{\mathcal {E}}}}
\def\tcF{{\widetilde{\mathcal {F}}}}
\def\tcG{{\widetilde{\mathcal {G}}}}
\def\tcH{{\widetilde{\mathcal {H}}}}
\def\tcI{{\widetilde{\mathcal {I}}}}
\def\tcJ{{\widetilde{\mathcal {J}}}}
\def\tcK{{\widetilde{\mathcal {K}}}}
\def\tcL{{\widetilde{\mathcal {L}}}}
\def\tcM{{\widetilde{\mathcal {M}}}}
\def\tcN{{\widetilde{\mathcal {N}}}}
\def\tcO{{\widetilde{\mathcal {O}}}}
\def\tcP{{\widetilde{\mathcal {P}}}}
\def\tcQ{{\widetilde{\mathcal {Q}}}}
\def\tcR{{\widetilde{\mathcal {R}}}}
\def\tcS{{\widetilde{\mathcal {S}}}}
\def\tcT{{\widetilde{\mathcal {T}}}}
\def\tcU{{\widetilde{\mathcal {U}}}}
\def\tcV{{\widetilde{\mathcal {V}}}}
\def\tcW{{\widetilde{\mathcal {W}}}}
\def\tcX{{\widetilde{\mathcal {X}}}}
\def\tcY{{\widetilde{\mathcal {Y}}}}
\def\tcZ{{\widetilde{\mathcal {Z}}}}

\def\tfA{{\widetilde{\mathfrak {A}}}}
\def\tfB{{\widetilde{\mathfrak {B}}}}
\def\tfC{{\widetilde{\mathfrak {C}}}}
\def\tfD{{\widetilde{\mathfrak {D}}}}
\def\tfE{{\widetilde{\mathfrak {E}}}}
\def\tfF{{\widetilde{\mathfrak {F}}}}
\def\tfG{{\widetilde{\mathfrak {G}}}}
\def\tfH{{\widetilde{\mathfrak {H}}}}
\def\tfI{{\widetilde{\mathfrak {I}}}}
\def\tfJ{{\widetilde{\mathfrak {J}}}}
\def\tfK{{\widetilde{\mathfrak {K}}}}
\def\tfL{{\widetilde{\mathfrak {L}}}}
\def\tfM{{\widetilde{\mathfrak {M}}}}
\def\tfN{{\widetilde{\mathfrak {N}}}}
\def\tfO{{\widetilde{\mathfrak {O}}}}
\def\tfP{{\widetilde{\mathfrak {P}}}}
\def\tfQ{{\widetilde{\mathfrak {Q}}}}
\def\tfR{{\widetilde{\mathfrak {R}}}}
\def\tfS{{\widetilde{\mathfrak {S}}}}
\def\tfT{{\widetilde{\mathfrak {T}}}}
\def\tfU{{\widetilde{\mathfrak {U}}}}
\def\tfV{{\widetilde{\mathfrak {V}}}}
\def\tfW{{\widetilde{\mathfrak {W}}}}
\def\tfX{{\widetilde{\mathfrak {X}}}}
\def\tfY{{\widetilde{\mathfrak {Y}}}}
\def\tfZ{{\widetilde{\mathfrak {Z}}}}

\def\Atil{{\widetilde A}}          \def\atil{{\tilde a}}
\def\Btil{{\widetilde B}}          \def\btil{{\tilde b}}
\def\Ctil{{\widetilde C}}          \def\ctil{{\tilde c}}
\def\Dtil{{\widetilde D}}          \def\dtil{{\tilde d}}
\def\Etil{{\widetilde E}}          \def\etil{{\tilde e}}
\def\Ftil{{\widetilde F}}          \def\ftil{{\tilde f}}
\def\Gtil{{\widetilde G}}          \def\gtil{{\tilde g}}
\def\Htil{{\widetilde H}}          \def\htil{{\tilde h}}
\def\Itil{{\widetilde I}}          \def\itil{{\tilde i}}
\def\Jtil{{\widetilde J}}          \def\jtil{{\tilde j}}
\def\Ktil{{\widetilde K}}          \def\ktil{{\tilde k}}
\def\Ltil{{\widetilde L}}          \def\ltil{{\tilde l}}
\def\Mtil{{\widetilde M}}          \def\mtil{{\tilde m}}
\def\Ntil{{\widetilde N}}          \def\ntil{{\tilde n}}
\def\Otil{{\widetilde O}}          \def\otil{{\tilde o}}
\def\Ptil{{\widetilde P}}          \def\ptil{{\tilde p}}
\def\Qtil{{\widetilde Q}}          \def\qtil{{\tilde q}}
\def\Rtil{{\widetilde R}}          \def\rtil{{\tilde r}}
\def\Stil{{\widetilde S}}          \def\stil{{\tilde s}}
\def\Ttil{{\widetilde T}}          \def\ttil{{\tilde t}}
\def\Util{{\widetilde U}}          \def\util{{\tilde u}}
\def\Vtil{{\widetilde V}}          \def\vtil{{\tilde v}}
\def\Wtil{{\widetilde W}}          \def\wtil{{\tilde w}}
\def\Xtil{{\widetilde X}}          \def\xtil{{\tilde x}}
\def\Ytil{{\widetilde Y}}          \def\ytil{{\tilde y}}
\def\Ztil{{\widetilde Z}}          \def\ztil{{\tilde z}}


\def\ahat{{\hat a}}          \def\Ahat{{\widehat A}}
\def\bhat{{\hat b}}          \def\Bhat{{\widehat B}}
\def\chat{{\hat c}}          \def\Chat{{\widehat C}}
\def\dhat{{\hat d}}          \def\Dhat{{\widehat D}}
\def\ehat{{\hat e}}          \def\Ehat{{\widehat E}}
\def\fhat{{\hat f}}          \def\Fhat{{\widehat F}}

\def\ghat{{\hat g}}          \def\Ghat{{\widehat G}}
\def\hhat{{\hat h}}          \def\Hhat{{\widehat H}}
\def\ihat{{\hat i}}          \def\Ihat{{\widehat I}}
\def\jhat{{\hat j}}          \def\Jhat{{\widehat J}}
\def\khat{{\hat k}}          \def\Khat{{\widehat K}}
\def\lhat{{\hat l}}          \def\Lhat{{\widehat L}}
\def\mhat{{\hat m}}          \def\Mhat{{\widehat M}}
\def\nhat{{\hat n}}          \def\Nhat{{\widehat N}}
\def\ohat{{\hat o}}          \def\Ohat{{\widehat O}}
\def\phat{{\hat p}}          \def\Phat{{\widehat P}}
\def\qhat{{\hat q}}          \def\Qhat{{\widehat Q}}
\def\rhat{{\hat r}}          \def\Rhat{{\widehat R}}
\def\shat{{\hat s}}          \def\Shat{{\widehat S}}
\def\that{{\hat t}}          \def\That{{\widehat T}}
\def\uhat{{\hat u}}          \def\Uhat{{\widehat U}}
\def\vhat{{\hat v}}          \def\Vhat{{\widehat V}}
\def\what{{\hat w}}          \def\What{{\widehat W}}
\def\xhat{{\hat x}}          \def\Xhat{{\widehat X}}
\def\yhat{{\hat y}}          \def\Yhat{{\widehat Y}}
\def\zhat{{\hat z}}          \def\Zhat{{\widehat Z}}


\def\intprod{\mathbin{\lr54}}
\def\reals{{\mathbb R}}
\def\integers{{\mathbb Z}}
\def\complex{{\mathbb C}\/}
\def\naturals{{\mathbb N}\/}
\def\distance{\operatorname{distance}\,}
\def\degree{\operatorname{degree}\,}
\def\dim{\operatorname{dimension}\,}
\def\Span{\operatorname{span}\,}
\def\ZZ{ {\mathbb Z} }
\def\e{\varepsilon}
\def\p{\partial}
\def\rp{{ ^{-1} }}
\def\Re{\operatorname{Re\,} }
\def\Im{\operatorname{Im\,} }
\def\ov{\overline}
\def\bx{{\bf{x}}}
\def\lt{L^2}
\def\Farrow{F} 
\def\Phiarrow{\Phi} 
\def\expect{\mathbb E}

\def\scriptx{{\mathcal X}}
\def\scriptb{{\mathcal B}}
\def\scripta{{\mathcal A}}
\def\scriptk{{\mathcal K}}
\def\scriptd{{\mathcal D}}
\def\scriptp{{\mathcal P}}
\def\scriptl{{\mathcal L}}
\def\scriptv{{\mathcal V}}
\def\scripti{{\mathcal I}}
\def\scripth{{\mathcal H}}
\def\scriptm{{\mathcal M}}
\def\scripte{{\mathcal E}}
\def\scriptt{{\mathcal T}}
\def\scriptb{{\mathcal B}}
\def\frakg{{\mathfrak g}}
\def\frakG{{\mathfrak G}}

\author{Michael Christ}
    \address{
        Michael Christ\\
        Department of Mathematics\\
        University of California \\
        Berkeley, CA 94720-3840, USA}
\email{mchrist@math.berkeley.edu}
\thanks{The first author was supported in part by NSF grant DMS-0401260}

\author{Andreas Seeger}
\address{
        Andreas Seeger\\
        Department of Mathematics\\
        University of Wisconsin\\
    Madison, Wisconsin 53706-1388, USA}
\email{seeger@math.wisc.edu}
\thanks{The second author was supported in part by NSF grant DMS-0200186}

\date{1-4-05}

\title[
Necessary conditions for vector-valued inequalities]
{
Necessary conditions\\ for  vector-valued operator inequalities\\ in
harmonic analysis}

\begin{abstract}
Via a  random construction we establish necessary conditions for
$L^p(\ell^q)$ inequalities for certain families of operators
arising in harmonic analysis. In
particular we consider dilates of a convolution kernel
with compactly supported Fourier transform,
vector maximal functions acting on classes of entire functions of
exponential type, and a characterization of Sobolev spaces by
square functions and pointwise moduli of smoothness.
\end{abstract}
\maketitle

\section{Introduction}
For $r>0$ let $\cE(r)$ be the space of all (smooth) distributions on
$\bbR^d$
whose Fourier transforms are supported in $\{\xi:|\xi|\le r\}$.
Also let $\eann(r)$   be the space of functions in
$\cE(r)$ whose Fourier transforms are supported in the annulus
$\{\xi:r/2\le |\xi|\le r\}$.

Let us first consider a convolution kernel $K$ whose Fourier
transform is compactly supported, say $K\in\cE(1)$. We are
concerned with vector valued inequalities involving dilates of
$K$, of the form
\begin{equation}\label{Lpellqineq}
\Big\|\Big(\sum_k|r_k^d K(r_k\cdot) *f_k|^q\Big)^{1/q}\Big\|_p\le
A \Big\|\Big(\sum_k|f_k|^q\Big)^{1/q}\Big\|_p.
\end{equation}
 An immediate necessary, but not sufficient, condition for
\eqref{Lpellqineq} to hold is that $K\in L^s$ for all $s\ge p$. This 
is seen by setting all but one $f_k$ to $0$ and 
(after  possibly a  rescaling) convolving
$K$ with a Schwartz function
 whose Fourier transform is equal to $1$ on the support of $\widehat K$.
In  the case $p>q$ we get a further necessary condition:

\begin{theorem}\label{thmconvol}
Suppose $0<q\le p<\infty$ and let
 $\{r_k\}_{k=1}^\infty$ be a fixed sequence of positive numbers.
Suppose that $K\in \cE(1)$ and
that \eqref{Lpellqineq} holds
 for all choices of  $f_k\in \cE(2r_k)$ with $\{f_k\}\in L^p(\ell^q)$.

Then $K\in L^q$  and
there exists a constant $C =C(p,q,d)$ so that
\begin{equation}\label{LqestimateforK}
\|K\|_q\le C(p,q,d) A;
\end{equation}
in particular $C(p,q,d)$ does not depend on the choice of the sequence
$\{r_k\}$.
\end{theorem}

As an application consider the Bochner-Riesz means defined by
$$\widehat {S_r^\la f}(\xi)= (1-r^{-2}|\xi|^2)^\la_+\widehat f(\xi).$$
Let $K_\la$ be the convolution kernel  for  $S^\la_1$.
From the well known formula for $K_\la$ (\cite{stein-weiss})
we know that $K_\la\in L^q$ if and only if $\la>d(1/q-1/2)-1/2$.
Consequently  if $q<p<2$ then
the operator
$$\{f_k\} \mapsto \{S_{r_k}^\la f_k\}$$
fails to be bounded on $L^p(\ell^q)$ if  $\la\le d(1/q-1/2)-1/2$, as well as on the  corresponding subspace with the restrictions
 $f_k\in \cE(2r_k)$.
This complements the familiar necessary condition
$\la>\max\{d(1/p-1/2)-1/2,0\}$
(\cite{stein-weiss}, \cite{fefferman}), which is known also to be sufficient for
certain $p$; for some refinements and implications
to known multiplier theorems see
 the remark at the end of \S\ref{sectionconvol} below.

We shall prove Theorem \ref{thmconvol}
  by a random construction which will be described in the
next section.
This construction  applies also to other situations, in particular  to
maximal functions which arise in the theory of function spaces. As the most
basic such example  we consider  a maximal operator acting on
functions of exponential type, which
was  introduced by Peetre \cite{Pe}, following  earlier related
research by Fefferman and Stein \cite{FS-Hp}.

For $r>0$ and  $\si\ge 0$ set
\begin{equation}\label{Peetre}
\fM_{\si,r} g(x)=\sup_y\frac{|g(x+y)|}{(1+r|y|)^\si}.
\end{equation}

As shown in \cite{Pe} one has the majorization
\begin{equation}\label{Msmajorization}
\fM_{\si,r} g(x)\lc [M_{HL}(|g|^s)]^{1/s}, \quad \forall \si\ge d/s,
\quad\text{ if } g\in \cE(r);
\end{equation}
here $M_{HL}$ denotes the Hardy-Littlewood maximal operator.
Now  by the Fefferman-Stein vector-valued maximal theorem (\cite{FS})
\begin{equation}
\label{vectorvaluedmaxfct}
\Big\|\Big(\sum_k|\fM_{\si,r_k} f_k|^q\Big)^{1/q}\Big\|_p
\lc
\Big\|\Big(\sum_k|f_k|^q\Big)^{1/q}\Big\|_p, \quad \si>
\max\{\frac dp, \frac dq\},
\end{equation}
provided that $f_k\in \cE(r_k)$ and
 $\{r_k\}_{k=1}^\infty$ is any   sequence of positive  radii.

It is well known
that the condition $\si>d/p$ is necessary --- again to see this
one simply chooses  a fixed Schwartz function for $g_1$ and sets $g_k=0$ for $k\ge 2$.
Moreover if $r_k=1$ for all $k$ the inequality clearly fails for all
$q\le p$; this is the same example that disproves an $L^{p}(\ell^1)$
inequality for the Hardy-Littlewood maximal function \cite{FS}. Indeed let
$\eta\in \cE(1)\cap \cS$, let $A$ be a large positive  integer,
and let $\{x(k)\}_{k=1}^{(2A+1)^d}$ be an enumeration of all integer lattice
points in the cube $Q_A$
of sidelength $A$ centered at the origin.
Define  $f_k(x)=\eta(x-x(k))$ if $1\le k\le (2A+1)^d$ and $f_k(x) =0$ otherwise.
Then $\|\{f_k\}\|_{L^p(\ell^q)}\lc A^{d/p}$.
Also  $\fM_{\si,1} f_k(x)\gc(1+|x-x(k)|)^{-\si}$ and a computation
shows that 
$\|\{\fM_{\si,1}f_k\}\|_{L^p(\ell^q)}\gc A^{d/p} \log A$ if $\si=d/q$ and
$\gc A^{d/p+d/q-\si}$ if $\si<d/q$.
Thus the condition
$\si>\max\{\frac dp, \frac dq\}$
in \eqref{vectorvaluedmaxfct}
 is sharp if $r_k\equiv 1$.

The preceding example does not immediately  apply to cases  where the
sequence of radii $r_k$ is sparse (say lacunary),
which happens in many of
the interesting cases for  which
\eqref{vectorvaluedmaxfct} is used.
Nevertheless we show that the condition
$\sigma>d/q$ is necessary for
\eqref{vectorvaluedmaxfct} to hold:

\begin{theorem}\label{thmpeetre}
Let  $\{r_k\}$ be any sequence of radii and suppose that $0< q\le p<\infty$.
Suppose $\si\le d/q$. Then there is a positive constant $c(p,q,\sigma,d)$ such that
for every $L\in \bbN$ there are functions
$f_k\in \cE_o(r_k)\cap L^p$, for $k=1,\dots, L$, so that
\begin{equation}\label{lowerboundmaxL}
\Big\|\Big(\sum_{k=1}^L|\fM_{\si,r_k} f_k|^q\Big)^{1/q}\Big\|_p
\ge c(p,q,\sigma,d) \max \{ L^{-\sigma+d/q} , \log^{1/q}\! L\}
\Big\|\Big(\sum_{k=1}^L|f_k|^q\Big)^{1/q}\Big\|_p.
\end{equation}
\end{theorem}
\noindent
Note that this lower bound holds for functions in $\cE_o(r_k)$, not merely in
$\cE(r_k)$.

\medskip

Next we shall state a  result on
a characterization of Sobolev spaces (or more general
Triebel-Lizorkin spaces)
by means of  pointwise moduli of continuity.
For $h\in \bbR^d$ let $\Delta_h f(x)= f(x+h)-f(x)$ and define higher difference operators inductively by $\Delta^0 f=f$,
$\Delta_h^m f=\Delta_h(\Delta_h^{m-1} f)$,
$m\ge 1$.
For suitable  classes of functions  let
\begin{equation}\label{definemaxderivative}
\fD_m^{\si,q} f(x)=
\Big(\int_0^1\frac{\sup_{|h|\le t} |\Delta_h^m f(x)|^q}{t^{1+\si q}} dt\Big)^{1/q}.
\end{equation}
It is known that if $m>\si$, $q=2$, and $1<p<\infty$ one can
characterize Sobolev spaces $\cL^p_\si$  using $\fD_m^{\si,2} $, namely
$\|f\|_{\cL^p_\si}:= \big\|\cF^{-1}[(1+|\cdot|^2)^{\si/2}\widehat f]\big\|_p \approx
\|f\|_p+ \|\fD_m^{\si,2} f\|_p
$
provided that $\si>\max\{d/p, d/2\}$.
This is a special case of a result
on Triebel-Lizorkin spaces $\fpqs$ (\cite{Tr1}, \cite{Tr2}).
We recall that $\fpqs$ is defined by
dyadic frequency decompositions; namely if $\beta_0\in \cE(1)$ so that
$\widehat \beta_0$ is equal to $1$ in a neighborhood of the origin, and if
$\beta_k= 2^{kd}\beta_0(2^k\cdot)-
2^{(k-1)d}\beta_0(2^{k-1}\cdot)$ for $k\ge 1$ then
$$
\|f\|_{\fpqs}\approx \Big\|\Big(\sum_{k=0}^\infty
2^{k\sigma q}|\beta_k*f|^q\Big)^{1/q}\Big\|_p;
$$
thus $F^{p2}_\si=\cL^p_\si$, $1<p<\infty$, by the usual  Littlewood-Paley
inequalities. Now by \cite{Tr1}, \S 2.5.10 we have for  $m>\si$,
$0<p<\infty$, $0<q\le \infty$ and  $\si> \max\{d/p, d/q\}$
\begin{equation}\label{fpqinequality}
\|f\|_p+ \|\fD_m^{\si,q} f\|_p \approx \|f\|_{\fpqs}\,.
\end{equation}


Again the condition $\si>d/p$ is necessary in
\eqref{fpqinequality},
 but it was  apparently open whether  for $p>q$
the characterization \eqref{fpqinequality} could hold without the
additional restriction $\si>d/q$ ({\it cf.} \cite{Tr2}). This was
pointed out to the second author by Herbert Koch  and Winfried
Sickel
 at an Oberwolfach  meeting some years ago.
We show that the restriction $\sigma>d/q$ is indeed necessary and
in the range $d/p<\sigma\le d/q$ we quantify the failure
of \eqref{fpqinequality} in terms of the support of the Fourier transform.

\begin{theorem}\label{thmsobolev}
Suppose that $0< \si<m $ and
 $0<q< p<\infty$.  For $r\ge 100$ let
$$\cA_{p,q,\sigma}(r)= \sup\big\{
\|\fD_m^{\si,q} f\|_p:\, \|f\|_{\fpqs}\le 1, f\in \cE(r) \big\}.
$$
Then
\begin{equation}
\label{equivalence-r}
\cA_{p,q,\si}(r)\approx (\log r)^{\frac dq -\sigma}\quad \text{ if } d/p<\sigma<d/q,
\end{equation}
and,
 for $\sigma=d/q$,
\begin{equation}
\label{equivalence-r-2}
\cA_{p,q,d/q}(r)\approx (\log\log r)^{1/q}
\quad \text{ if }  q\le 1.
\end{equation}
Moreover,
\begin{equation}
\label{inequality-r}
C^{-1} (\log\log r)^{1/q}\le
\cA_{p,q,d/q}(r) \le C \log\log r \quad \text{ if } 1<q<p.
\end{equation}
\end{theorem}

In \eqref{equivalence-r} the notation $a_1\approx a_2$ means that
there is a positive constant $C=C(p,q,d,\sigma,m)$ which does not
depend on $r$
so that $C^{-1}a_1\le a_2\le Ca_1$.
An application of the Banach-Steinhaus theorem
(\cf. Theorems 2.5, 2.6
in \cite{Rudin}) shows that
for  $\si\le d/q$ there is
an $f\in \fpqs(\bbR^d)$ for which
$\fD_m^{\si,q} f$ does not belong to $L^p(\bbR^d)$
(in fact this holds for a class of second category in $\fpqs$).


Finally we settle an endpoint question about oscillatory multipliers on
the $F$-spaces. Consider  the operator given by
\begin{equation}
\label{oscillatorymult}
\widehat {T_{\gamma,b} f}(\xi)=
\frac{e^{i|\xi|^{\gamma}}}{(1+|\xi|^2)^{b/2}}\widehat f(\xi),
\end{equation}
for $0<\gamma<1$.
It is well known that
$T_{\gamma,b}$ maps the Besov spaces $B^p_{\alpha,q}$ into itself
if and only $b/\gamma\ge   d|1/p-1/2|$,
and by a simple application of H\"older's inequality
the same result holds for $F^p_{\alpha,q}$ with the strict inequality
$b/\gamma>d|1/p-1/2|$. If
 $1\le p\le q\le p'$ (if $p>1$), or  $p\le q\le \infty$, $p<1$ the  endpoint result with $b/\gamma=d|1/p-1/2|$  holds for the $F$-spaces
(see \cite{FS-Hp}, and for more general multiplier theorems
\cite{BS}, \cite{SeStud}).

We show that for the endpoint result the restriction on $q$ is necessary.

\begin{theorem}\label{oscillatorytheorem}
Let  $0<q<p\le 2$, $\alpha\in \bbR$ and let  $0<\gamma<1$, $b=\gamma d(1/p-1/2)$.
Then for $r\ge 2$
\begin{equation}
\label{oscequiv}
\sup
\big\{\|T_{\gamma,b}f\|_{F^{p}_{\alpha,q}}: \|f\|_{F^{p}_{\alpha,q}}\le 1, f\in \cE(r)
\big\}
\approx \big(\log r\big)^{1/q-1/p}
\end{equation}
\end{theorem}

\medskip

In \S\ref{sectionrandom} we shall give the basic random
construction that underlies the proofs of all the theorems.
Theorem \ref{thmconvol}
is proved in \S\ref{sectionconvol}. Theorem \ref{thmpeetre} will be proved in
\S\ref{sectionpeetre} and a second deterministic proof of  the lacunary case
will be given in  \S\ref{deterministic}.
Theorem \ref{thmsobolev}  will be proved in \S\ref{pfthmsobolev}
and Theorem \ref{oscillatorytheorem} in \S\ref{pfoscillatory}.

\section{A random construction}\label{sectionrandom}

For each $n\in\{0,1,2,\cdots\}$
let $\cQ(n)$ be the set of all dyadic cubes of sidelength $2^{-n}$
in $[0,1)^d$; more specifically
all cubes of the form
$\prod_{i=1}^d[j_i 2^{-n},(j_i+1) 2^{-n})\,$
where the $j_i$ are integers, $0\le j_i<2^{n}$, for $i=1,\dots, d$.
For any dyadic cube $Q$ let  $\chi_{Q}$
denote the characteristic function of $Q$.

Let $a\in (0,1)$ be a parameter to be specified.
Let $\Omega$ be  a probability space  with probability measure $\mu$,
on which there is a family $\{\theta_{Q,a}\}$ of independent
random variables indexed by the dyadic subcubes of $[0,1]^d$,
each of which takes the value $1$ with probability $a$ and the value $0$
with probability $1-a$.
If $B\subset\Omega$ we denote by $\mu(B)$ the probability
of $B$ and the  expectation of a function $g$ on $\Omega$
(i.e.\ a random variable) is  given by the integral
$\expect(g) = \int_\Omega g(\omega)\,d\mu(\omega)$.

In what follows we fix a sequence
$\{n_k\}_{k=0}^\infty$ of nonnegative integers.
We consider random functions
\begin{equation}\label{hk}
h_{k}^{\om,a}(x)
=
\sum_{Q\in \cQ(n_k)}
\theta_{Q,a}(\omega)\chi_{Q}(x);
\end{equation}
these are supported on $[0,1]^d$. Note that $h_{k}^{\om,a}(x)\in \{0,1\}$
 for all $x$.
The parameter $a$ will be mostly fixed (except in \S7), and we use the notation
$h^\om_k\equiv h^{\om,a}_{k}$, $\theta_{Q}=\theta_{Q,a}$ if
the value of $a$ is clear.

\begin{lemma}\label{bernoulliupperestimate}
Suppose $p,q\in (0,\infty)$ and   $0<a< C_1 L^{-1}$. Let $\sigma>
\max\{d/p, d/q\}$.
Then
\begin{equation} \label{Homupperbound}
\Big(\int_\Om \Big\| \Big(\sum_{k=1}^L
[\fM_{\sigma,2^{n_k}}h^{\om,a}_{k}]^q\Big)^{1/q}\Big\|_p^p
d\mu\Big)^{1/q}
\le C(p,q,C_1)
\end{equation}
\end{lemma}

\begin{proof}
We first observe that  for $r>0$ and
every $x\in [0,1]^d$,
\begin{equation}\label{pointwisebernoulli}
\int_\Omega \Big(\sum_{k=1}^L |h_{k}^\om(x)|^q\Big)^{r} d\mu   =
\sum_{n=0}^{L} \binom{L}{n} a^n(1-a)^{L-n} n^{r}.
\end{equation}
To see this let  $x\in [0,1)^d$ and observe that
for each $k$,
$h_{k}^\om(x) = \theta_{Q}(\omega)\chi_{Q}(x)$
for a single  $Q=Q(k,x)\in \cQ(n_k)$ and thus also $h_{k}^\om(x)
=[h_{k}^\om(x)]^q$.
One has then $2^L$ possible events, indexed by all subsets
$S\subset\{1,2,\cdots,L\}$;
the event $\Omega(S,x)$  that $\theta_{Q(k,x)}(\omega)$ equals $1$
 for all $k\in S$  and equals $0$  for all $k\notin S$
has probability
$$\mu (\Omega(S,x))=a^{\card(S)}(1-a)^{L-\card(S)},$$ by independence.
The function $(\sum_{k=1}^L h_k(x,\omega))^r$ has value $\card(S)^r$ at
such an event.
Lastly the number of subsets $S$ having cardinality $n$ is $\binom{L}{n}$.
Thus, for
every $x$,
\begin{align}
\int_\Omega\Big( \sum_{k=1}^L h_k^\om(x)\Big)^{r}d\mu
&=
\sum_{n=0}^{L} \sum_{\card(S)=n}
\int_{\Omega(S,x)}\Big(\sum_{k\in S} h_k^\om(x)\Big)^{r}d\mu
\notag
\\&=
\sum_{n=0}^{L} \binom{L}{n} a^n(1-a)^{L-n} n^r
\label{binomialsum}
\end{align}
which  gives
\eqref{pointwisebernoulli}.

We set $r=p/q$ in \eqref{pointwisebernoulli} and let
  $r_0$ be the smallest positive
integer $\ge p/q$.
Then
\begin{equation*}
\sum_{n=0}^{L} \binom{L}{n}a^n (1-a)^{L-n} n^{p/q}
\le
\sum_{n=1}^{\infty} \frac{L^n}{n!} a^n  n^{r_0}
=
(t\frac{d}{dt})^{r_0} e^t\Big|_{t=La} \le C(r_0,C_1).
\end{equation*}
By
\eqref{pointwisebernoulli}, integration in $x$ and Fubini's theorem the last inequality implies
\begin{equation}\label{expectedbound}
\Big(\int_\Omega \Big\|\Big(\sum_{k=1}^L |h_k^\om|^q\Big)^{1/q}\Big\|_p^p
 d\mu \Big)^{1/p}  \le C(p/q, C_1)\quad  \text{ if } La\le C_1.
\end{equation}

The conclusion of the lemma now follows from \eqref{vectorvaluedmaxfct},
but we repeat the derivation since it involves an estimate that will be needed later.
Observe that since $h_k^\om$ assumes only the values $1$ and $0$ and is constant on dyadic cubes of length $2^{-n_k}$
there is the estimate
\begin{equation}\label{msestimate}
\sup_{2^{-n_k+l}\le |y|\le 2^{-n_k+l+1}}|h_{k}^\om(x+y)|\le C_s
 2^{l \frac ds} \big(M([h_{k}^\om]^s)\big)^{1/s}
\end{equation}
 for any $s\le 1$. Consequently
$\fM_{\sigma,r_k}
[h_{k}^\om](x) \le C_\sigma  \big(M([h_{k}^\om]^s)\big)^{1/s}$
if $\sigma>d/s$ and the
vector Fefferman-Stein inequality \cite{FS} can be applied if
$p/s>1$, $q/s>1$. Thus
the asserted maximal inequality follows from \eqref{expectedbound}.
\end{proof}

An immediate consequence is

\begin{corollary}\label{bernoulliupperestimateeta}
Suppose $p,q\in (0,\infty)$,
$\sigma>\max\{d/p, d/q\}$,
$L\in\bbN$ and   $0<a< C_1 L^{-1}$.

Let $\eta$ be a Schwartz function and
$\eta_k(x)=2^{n_k d}\eta(2^{n_k}x).$
Denote by
$\Farrow^{\omega}$
  the random vector-valued
function  defined by
$\Farrow^{\omega}_{k}(x) =  \eta_k*h_{k}^{\omega,a}(x)$  if  $ 1\le k\le L$,
and
$\Farrow^{\omega}_{k}(x) =0$  if $ k>L$.
Then
\begin{equation} \label{Fomupperbound}
\Big(\int_\Om \big\|F^{\om}\|_{L^p(\ell^q)}^p d\mu\Big)^{1/p}
\le \widetilde C(p,q,C_1)
\end{equation}
\end{corollary}

\medskip

\noi{\it Remark.}
The quantity
\eqref{binomialsum} is bounded by
$ C_r La$ if $L^{-1}\le a\le 1$, see a calculation in Bourgain \cite{bourgain}. There is also a corresponding lower bound  for $r\ge 1$, in fact
there is the identity
$\sum_{n=0}^{L} \binom{L}{n} b^n(1-b)^{L-n} n = Lb, \quad 0<b<1.$
To see this observe that the left hand side
is equal to
$(1-b)^L t\tfrac{d}{dt}(1+t)^L$ when evaluated at $t=b/(1-b)$.
One also has
$(\sum_{n=0}^{L} \binom{L}{n} b^n(1-b)^{L-n} n^r)^{1/r}\ge Lb$ if $r\ge 1
$; this follows from H\"older's inequality since
$\sum_{n=0}^{L} \binom{L}{n} b^n(1-b)^{L-n}=1$.

\section{Proof of Theorem \ref{thmconvol}}
\label{sectionconvol}

For $z\in \bbR^d$, $\ell\in \bbZ$ denote by $\fQ(\ell,z)$
 the family of all cubes of the form $z+Q$, with $Q$  any  dyadic cube
of sidelength $2^{-\ell}$ in $\bbR^d$.
We shall use the important Plancherel-P\'olya theorem for entire functions
of exponential type
(\cite{Pl-Po}, \cite{Tr1}).
It
says that there are absolute positive constants $C$, $m$ depending only on $q\in (0,\infty)$ and $d$  so that for all $\ell$, $z$
\begin{equation}\label{PP}
C^{-1} \Big(\sum_{Q\in \fQ(\ell+m,z)}|f(x_Q)|^q\Big)^{1/q}
\le 2^{\ell d/q}\|f\|_q \le
C \Big(\sum_{Q\in \fQ(\ell+m,z)}|f(\widetilde x_Q)|^q\Big)^{1/q},
\quad f\in  \cE(2^\ell);
\end{equation}
here $x_Q\in Q$, $\widetilde x_Q\in Q$ and the constants in \eqref{PP} are independent of the specific choices of $x_Q$,
$\widetilde   x_Q$.

An equivalent formulation is
\begin{equation}\label{PPalt}
C^{-1} \Big(\int\sup_{|x-y|\le u 2^{-k}} |f(y)|^q dx\Big)^{1/q}
\le \|f\|_q \le
C \Big(\int\inf_{|x-y|\le u 2^{-k}} |f(y)|^q dx\Big)^{1/q},
\quad f\in  \cE(2^k);
\end{equation}
here $C$ and  $u\in (0,1)$
 depend only on $q$ and $d$.

As the statement of Theorem \ref{thmconvol} is trivial for $p\le q$ we
shall assume $p\ge q$ in what follows.
If  $K\in \cE(1)$ satisfies condition \eqref{Lpellqineq} with
$q\le p$
we shall show that for all $N\in \bbN$
\begin{equation}\label{discretelowerbound}
\Big(\int_{|x|\le 2^N}\inf_{|x-y|\le u } |K(y)|^q dx\Big)^{1/q}
\le C(q,d,u) A.
\end{equation}
Here  we may pass to the limit as $N\to \infty$ and then,
  choosing $u=u(q,d)$
sufficiently small,  we may apply the  second inequality in \eqref{PPalt}
 to deduce the assertion of Theorem \ref{thmconvol}. In what follows we pick an integer
$M$ so that $2^{-M+d+1}\le u<2^{-M+d+2}$.

In order to show \eqref{discretelowerbound}
we may  use
\eqref{Lpellqineq}
 for functions $\{f_k\}_{k=1}^L$ indexed by a finite family
of radii; we put $L=2^{Nd}$ and  by a scaling we may assume
that \begin{equation} \label{rkcondition}
r_k\ge 2^{10 d+10 N}, \qquad k=1,\dots, L.
\end{equation}
It will be useful to replace $K$ with a kernel which vanishes for
$|x|\ge 2^{N+2}$.
Let $\zeta$ be a $C^\infty$ function with
compact support in $\{x:|x|<4\}$ which
equals $1$ for $|x|\le 2$. Let
$\zeta_N(x)=\zeta(2^{-N}x)$ and let $K^N=K\zeta_N$.
Clearly
\eqref{discretelowerbound} follows from
\begin{equation}\label{discretelowerboundKN}
\Big(\int \inf_{|x-y|\le u } |K^N(y)|^q dx\Big)^{1/q}
\le C'(q,d,u) A.
\end{equation}

We first deduce from \eqref{Lpellqineq}
a vector-valued inequality
for the dilates of $K^N$. We define positive integers $n_k$ as in the
 previous section,
namely by
\begin{equation}
2^{n_k-M-d-1}\le r_k<2^{n_k-M-d}.
\end{equation}
With these specifications on $r_k$, $n_k$ we prove

\begin{lemma}\label{LqineqKNlemma}
 Suppose that $q\le p$ and   that \eqref{Lpellqineq} and
\eqref{rkcondition} hold. Set $K_k^N(x)=r_k^dK^N(r_k x)$,
$a=L^{-1}=2^{-Nd}$ and define $h_k^\om\equiv h_{k}^{\om,a}$ as in
\eqref{hk}.
Then
\begin{equation}\label{LqineqKN}
\Big(\sum_{k=1}^L\int_{[0,1]^d}\int_\Omega |K_k^N*h_{k}^{\om,a}|^q d\mu \,
dx\Big)^{1/q} \le CA
\end{equation}
\end{lemma}

\begin{proof} By H\"older's inequality and Fubini's theorem
\begin{equation}\label{HolderKN}
\Big(\sum_{k=1}^L\int_{[0,1]^d}\int_\Omega |K_k^N*h_{k}^{\om}|^q d\mu
dx\Big)^{1/q} \le
\Big(\int_\Omega \int_{[0,1]^d}\Big(\sum_{k=1}^L|K_k^N*h_{k}^{\om}|^q
\Big)^{p/q} dx d\mu  \Big)^{1/p}.
\end{equation}

Let
$e_z(x)=e^{i\inn xz}$. Then for any compactly supported bounded function $g$
\begin{equation} \label{averaging}
 K^N_k* g(x) = (2\pi)^{-d} \int \widehat {\zeta_N}(\xi)
 e_{r_k\xi}(x)
K_k*   [g e_{-r_k\xi}] (x)d\xi.
\end{equation}
Let $\eta$ be a Schwartz function in $\cE(2)$ with the property that
$\widehat \eta(\xi)=1$ for $|\xi|\le 1$. Let $\eta_k=r_k^d\eta(r_k\cdot)$
and
 $K_k=r_k^dK(r_k \cdot)$,  then
\begin{equation}\label{reproducing}
K_k*\eta_k=K_k.
\end{equation}

Now suppose  $1\le q\le p$. Then \eqref{averaging}, \eqref{reproducing},
 Minkowski's inequality and  the assumption
\eqref{Lpellqineq} imply for fixed $\om$
\begin{align*}
\Big\| \Big(\sum_{k=1}^{L}| K_k^N* h_{k}^\om|^q\Big)^{1/q}\Big\|_p
&\le \int|\widehat {\zeta_N}(\xi)| \Big\| \Big(\sum_{k=1}^{L}\big|
K_k* \eta_k*[h_{k}^\om e_{-r_k \xi}]\big|^q\Big)^{1/q}\Big\|_p d\xi
\\
&\le A\int|\widehat {\zeta_N}(\xi)| \Big\| \Big
(\sum_{k=1}^{L}\big|\eta_k*[h_{k}^\om e_{-r_k \xi}]\big|^q\Big)^
{1/q}\Big\|_p d\xi
\\
&\le C_\rho A
\Big\| \Big(\sum_{k=1}^{L}
 \Big[\sup_y\frac{|h_{k}^\om(\cdot +y)|}{(1+r_k|y|)^\rho}
\Big]^q\Big)^{1/q}\Big\|_p
\end{align*}
for any $\rho>0$. We have used that $\|\widehat {\zeta_N}\|_1=O(1)$.
We choose  $\rho>d/q$, take $p$th powers, and integrate over
$\omega\in \Omega$. By Lemma \ref{bernoulliupperestimate} we obtain
\begin{equation} \label{LpellqKN}
\Big(\int_\Omega\Big\|
\Big(\sum_{k=1}^{L}| K_k^N* h_{k}^\om|^q\Big)^{1/q}\Big\|_p^pd\mu
\Big)^{1/p}\le C(p,q,d) A
\end{equation}
and \eqref{LqineqKN}
follows from \eqref{LpellqKN}
and \eqref{HolderKN} (in the case $q\ge 1$).

It remains to prove \eqref{LpellqKN} in the case
$q\le 1$. Since $\zeta$ has compact support we can apply the
Plancherel-P\'olya theorem in $L^q$.
Let $\{Q_\nu^N\}$ denote the collection of dyadic cubes of sidelength
$2^{-M-N}$ (where $2^{-M}\approx u=u(q)$ as in \eqref{PPalt}). For each
 such cube choose
$\xi_\nu\in Q_\nu^N$. Then for fixed $\om$
\begin{align*}
\Big\| \Big(\sum_{k=1}^{L}\big| K^N_k* h_k^\om\big|^q\Big)^{1/q}\Big\|_p
&\le \Big(\int \Big(\sum_{k=1}^{L} \Big(\int\big|\widehat
{\zeta_N}(\xi)K_k* [h_{k}^\om e_{-r_k\xi}](x)\big
|d\xi\Big)^q\Big)^{p/q}dx\Big)^{1/p}
\\
&\lc \Big(\int \Big(\sum_{k=1}^{L} \Big(\sum_\nu 2^{-Nd}|\widehat
{\zeta_N}(\xi_\nu)K_k*[h_{k}^\om e_{-r_k\xi_\nu}](x)|\Big)^q\Big)^{p/q}dx\Big)^{1/p}
\end{align*}
and by the imbedding $\ell^q\subset \ell^1$ and
 Minkowski's inequality ($p/q\ge 1$)
 this is dominated by
\begin{align*} &\Big(\int \Big(\sum_{k=1}^{L} \sum_\nu 2^{-Ndq}\big|\widehat
{\zeta_N}(\xi_\nu)
K_k*[h_{k}^\om e_{-r_k\xi_\nu}](x)\big|^q\Big)^{p/q}dx\Big)^{1/p}
\\
&\lc  \Big(2^{-Ndq} \sum_\nu |\widehat {\zeta_N}(\xi_\nu)|^q
\Big(\int\Big(\sum_{k=1}^L\big|
K_k*[h^\om_{k}e_{-r_k\xi_\nu}(x)]\big|^q\Big)^{p/q}dx\Big)^{q/p}\Big)^{1/q}.
\end{align*}
By \eqref{reproducing} and \eqref{Lpellqineq} the last expression is
 in turn  dominated by
\begin{align}
&\Big(2^{-Ndq} \sum_\nu |\widehat {\zeta_N}(\xi_\nu)|^q
A^q\Big(\int\Big(\sum_{k=1}^L\big|
\eta_k*[h^\om_{k}e_{-r_k\xi_\nu}]\big|^q\Big)^{p/q}dx\Big)^{q/p}\Big)^{1/q}
\notag
\\
&\lc C_M A \Big(2^{-Ndq}
\sum_\nu|\widehat {\zeta_N}(\xi_\nu)|^q\Big)^{1/q}
\Big\|\Big(\sum_{k=1}^L
\big |\sup_y\frac{|h_{k}^\om(\cdot +y)}{(1+r_k|y|)^\rho}
\big|^q\Big)^{1/q}\Big\|_p.
\label{lastboundinlemma}
\end{align}
To eliminate  the $\nu$-summation we observe that by
the Plancherel-P\'olya  theorem
$$
2^{-Ndq}
\sum_\nu|\widehat \zeta_N(\xi_\nu)|^q\lc
2^{-Nd(q-1)}
\int|\widehat \zeta_N(\xi)|^q d\xi =\int|\widehat\zeta(\xi)|^q  d\xi.
$$
Thus we may
 apply Lemma  \ref{bernoulliupperestimate} (choosing  $\rho>d/q$)
to bound \eqref{lastboundinlemma} and  obtain
\eqref{LpellqKN} in the case $q<1$ as well.
\end{proof}

\noi{\bf Proof of  Theorem \ref{thmconvol}, conclusion.}
Let $Q^{N+M+2}_k(x)$ be the unique dyadic cube of sidelength $2^{-n_k+N+M+2}$
containing $x$ and let
$V^{N,M}_k(x)$ be the union of all dyadic  cubes of sidelength
$2^{-n_k+N+M+2}$ whose boundaries have
nonempty intersection with the boundary of $Q^{N+M+2}_k(x)$.
Then $V^{N,M}_k(x) \subset [0,1]^d$ provided that $x\in [1/4,3/4]^d$.
Let $\cV^{N,M}_k(x)$ be the family of all dyadic  cubes in $\cQ(n_k)$
which are contained in the closure of $V^{N,M}_k(x)$.

One of the obstacles to be overcome in our proofs is that unwanted
cancellations could conceivably arise between the different terms
contributing to expressions such as $$\sum_{Q\in\cQ(n_k)}
\theta_Q(\omega) K^N_k*\chi_Q(x).$$ We will handle this by
considering the contributions of events in which all terms but one
in the sum are either small, or have coefficients
$\theta_Q(\omega)=0$. To this end, for each $Q\in\cV^{N,M}_k(x)$
define the event
\begin{equation}
\Om(k,x,Q)=\{\omega\in\Omega: \text{ $\theta_Q(\om)=1$ and
$\theta_{Q'}(\om)=0$ for all $Q'\in
\cV^{N,M}_k(x)\setminus\{Q\}$\} }.
\end{equation}

If $Q\in \cQ(n_k)$ but $Q\notin \cV^{N,M}_k(x)$ then
$r_k|x-y|\ge r_k2^{-n_k+N+M+2}\ge  2^{N+2}$
for all $y\in Q$ and thus $K^N_k*\chi_Q=0$.
For fixed $1\le k\le L$, $x\in [1/4, 3/4]^d$,
\begin{align*}
\int_\Omega |K_k^N*h_{k}^\om(x)|^q d\mu
&=\int_\Omega\Big|
\sum_{
 Q'\in \cV^{N,M}_k(x) }
\theta_{Q'}(\om) K^N_k*\chi_{Q'}(x)\Big|^q d\mu
\\
\ge &\sum_{ Q\in \cV^{N,M}_k(x) }
\int_{\Omega(k,x,Q)}\Big|
\sum_{Q'\in \cV^{N,M}_k(x) }
\theta_{Q'}(\om) K^N_k*\chi_{Q'}(x)\Big|^q d\mu
\\
= &\sum_{Q\in \cV^{N,M}_k(x) }
\mu(\Om(k,x,Q))\big|K^N_k*\chi_{Q}(x)\big|^q.
\end{align*}
Now  \begin{equation*} \mu(\Omega(k,x,Q))
=a(1-a)^{2^{(N+M+2)d}3^d -1}
\end{equation*}
 with $a=2^{-Nd}$. Thus
\begin{equation}\mu(\Omega(k,x,Q)) \ge c_M 2^{-Nd}\end{equation} and therefore
\begin{equation*}
\mu(\Om(k,x,Q))\big|K^N_k*\chi_{Q}(x)\big|^q
\ge c_M 2^{-Nd} 2^{-n_k dq}
\inf_{y\in Q}|r_k^dK^N(r_k(x-y))|^q.
\end{equation*}
We have thus proved that
\begin{equation*}
\int_\Omega |K_k^N*h_{k}^\om(x)|^q d\mu
\ge c_M' 2^{-Nd}
\sum_{Q\in \cV^{N,M}_k(x) }
\inf_{z\in r_kQ}|K^N(r_kx -z))|^q.
\end{equation*}
Now the disjoint  cubes $r_k x-r_kQ$
cover  the ball of radius $2^{N+2}$ as $Q$
ranges over the cubes in $\cV^{N,M}_k(x)$.
The diameter of $r_kx-r_kQ$ is bounded by $ \sqrt
d r_k 2^{-n_k}\le 2^{-M}\le u$
and therefore
\begin{equation}\label{fixedxkconclusion}
\int_\Omega |K_k^N*h_{k}^\om(x)|^q d\mu
\ge c_M' 2^{-Nd}
\int \inf_{|y-z|\le u}|K^N(z)|^q dy.
\end{equation}
Now integrate over $x\in [1/4/3/4]^d$ and sum in $k=1,\dots, 2^{Nd}$
and the assertion
\eqref{discretelowerboundKN} follows from
\eqref{fixedxkconclusion} and \eqref{LqineqKN}.\qed

\medskip

\noi{\it Remark.}
Theorem \ref{thmconvol} can be applied to the case  of Bochner-Riesz multipliers mentioned in the introduction. A refinement of this example is as  follows.
Let $\chi$ be supported in $\{\xi:3/4<|\xi|<5/4\}$
 and be equal to $1$  in a neighborhood of the unit circle  and
consider the multiplier $$m_{\la,\delta} (\xi)=
\sum_{k\in \bbZ} \chi(2^{-k}\xi)
(1-2^{-2k}|\xi|^2)^\la_+  [\log(1-2^{-k}|\xi|)^{-1} ]^{-\delta}$$
Then
$f \mapsto \cF^{-1}[m_{\la,\delta} \widehat f]$ fails to be bounded on
the homogeneous Triebel-Lizorkin
space
${\dot F}^{pq}_0$ if
$\la<d(1/q-1/2)-1/2$, or
$\la=d(1/q-1/2)-1/2$, $\delta\le 1/q$.
These  examples show that that  the
 restriction $p\le q\le p'$ in some multiplier theorems for
Triebel-Lizorkin spaces stated in \cite{SeTams}, \cite{SeStud} is
needed; moreover, if $q\le 1$ then the  condition on $q$ in the
analogue of the Mikhlin-H\"ormander multiplier theorem stated on
p.75 in \cite{Tr1} is necessary.

\section{Proof of Theorem  \ref{thmpeetre}}\label{sectionpeetre}

We use the random construction of \S\ref{sectionrandom}. Fix  a
real valued Schwartz function $\eta$ so that
$\widehat \eta$ is supported  in $\{\xi: 1/2<|\xi|<1\}$ and so that
 $\eta(x)\ge 1$ for $|x|\le 2^{-M+2+d}$ (with some positive $M$
which is fixed in the proof).

Let $\sigma\le d/q$ and $q\le p$ and let $L$ be large.
We
may assume that $L=2^{Nd}$ for some large $N\in \bbN$. To show the lower bound
\eqref{lowerboundmaxL} we may assume that the $r_k$'s, $k=1,\dots, L$ are
large.
This follows by scaling, namely if
 $\delta_tf(x):=f(tx)$, and  if $f_k\in \cE_o(r)$ then
 $\delta_t f_k\in \cE_o(tr)$; moreover
$\delta_t^{-1}\fM_{\sigma,rt}\delta_t =
\fM_{\sigma,r}$. Thus
the  operator  norms of $\{\fM_{\sigma,r_k}\}$ and
$\{\fM_{\sigma,t r_k}\}$ are the same.

We may assume
$$2^{n_k-M}\le r_k<2^{n_{k}+1-M}, \quad n_k\in \bbZ, \quad
n_k\ge 100 d +M+N,
$$
for  $k=1,\dots, L$.
Define
$\eta_k(x) =r_k^d \eta(r_k x)$ and
\begin{equation}\label{definitionofgkom}
g_{k}^{\om,a}= \eta_k*h_{k}^{\om,a}, \quad a=2^{-Nd},
\end{equation}
with $h_{k}^{\om,a}$ as in
\eqref{hk}. Note that
 $g_{k}^{\om,a}\in \cE_o(r_k)$. We  omit the superscript $a$ in what follows.

Since $p\ge q$ we see  by H\"older's inequality and Fubini's theorem that
\begin{equation}\label{holderineq}
\Big(\int_\Om\Big\|\Big(\sum_{k=1}^L
|\fM_{\si,r_k} g_{k}^\om|^q\Big)^{1/q}\Big\|_p^p d\mu\Big)^{1/p}\ge
\Big(\sum_{k=1}^L\int_{[0,1]^d}\int_\Om
|\fM_{\si,r_k} g_{k}^\om|^qd\mu \,dx\Big)^{1/q}.
\end{equation}

Let $x\in [1/4,3/4]^d$ and let $Q^j_k(x)$
be the unique dyadic cube of sidelength $2^{-n_k+j}$
which contains $x$.
Let $$\cM_{j,k} f(x)=\sup_{y\in  Q^j_k(x)\setminus Q^{j-1}_k(x)}|f(y)|$$
Then
\begin{equation}\label{Mjkcomparison}
\fM_{\sigma,r_k}[g_{k}^\om](x)\ge c_{M,d}\sup_{2\le j\le N}
2^{-j\sigma} \cM_{j,k}g_{k}^\om(x).
\end{equation}
Thus, in view of Corollary \ref{bernoulliupperestimateeta},
\eqref{holderineq} and  \eqref{Mjkcomparison}
 it suffices to show that for $\sigma\le d/q$
\begin{equation}\label{claimoflowerbound}
\Big(\sum_{k=1}^L\int_{[\tfrac 14,\tfrac 34]^d}\int_\Om
\sup_{2\le j\le N}|2^{-j\sigma} \cM_{j,k}g_{k}^\om(x)
|^q d\mu dx\Big)^{1/q}\ge c \max\{2^{N(-\sigma+d/q)}, N^{1/q}\}.
\end{equation}

To show \eqref{claimoflowerbound} we let
$V^N_k(x)$ be the union of all dyadic  cubes of sidelength
$2^{-n_k+N+1}$ whose boundaries have
nonempty intersection with the boundary of $Q^{N+1}_k(x)$.
Split
\begin{equation*}
\cM_{j,k}g_{k}^\om(x) =I^\om_{j}(k,x)+II^\om_j(k,x)
\end{equation*}
where
\begin{align*}
I^\om_j(k,x)&=
\sup_{y\in  Q^j_k(x)\setminus Q^{j-1}_k(x)}\Big|
\sum_{Q\in \cQ(n_k) \atop {Q\subset V^N_k(x)}}
\theta_Q(\om)\eta_k*\chi_Q(y)\Big|
\\
II^\om_j(k,x)&=
\sup_{y\in  Q^j_k(x)\setminus Q^{j-1}_k(x)}\Big|
\sum_{Q\in \cQ(n_k) \atop Q\subset [\tfrac 14,\tfrac34]^d\setminus  V^N_k(x)}
\theta_Q(\om)\eta_k*\chi_Q(y)\Big|
\end{align*}
The terms $II^\om_j(k,x)$ are error terms; indeed if $y\in Q^j_k(x)$, $j\le N$
and  $z\in[\tfrac 14,\tfrac34]^d\setminus  V^N_k(x)$ then $|y-z|\ge c2^{-n_k+N}$ and from this it is easy to see that
\begin{equation*}
\sup_{2\le j\le N}\big|II^\om_j(k,x)
|\le C_{M,\rho} \fM_{\rho, 2^{-n_k}}
[h_{k}^\om]
\end{equation*}
for any $\rho>0$. Thus by
Lemma \ref{bernoulliupperestimate}
\begin{equation}\label{IIest}
\Big(\sum_{k=1}^L\int_{[\tfrac 14,\tfrac 34]^d}\int_\Om
\sup_{2\le j\le N}|2^{-j\sigma}II^\om_j(k,x)|^q d\mu \,dx\Big)^{1/q}\le C(p,q,M).
\end{equation}

We show for almost every $x\in [1/4,3/4)^d$, $1\le k\le L$
the uniform lower bound
\begin{equation}\label{I-lowerbd}
\int_\Om
\sup_{2\le j\le N}|2^{-j\sigma} I^\om_j(k,x)
|^q  d\mu  \ge c' 2^{-Nd}
\max\{2^{N(d-q\sigma)},  N\}.
\end{equation}

Clearly \eqref{claimoflowerbound} follows from
\eqref{I-lowerbd} after integrating in $x$ and then
summing in $k$ (recall that $L=2^{Nd}$); the error term \eqref{IIest}
changes this lower bound only by a small constant if $N$ is large.

Next,  to prove \eqref{I-lowerbd} we observe that if
 $Q\in \cQ(n_k)$ and if $y_Q$ is the center of $Q$ and $z\in Q$
then $|y_Q-z|\le \sqrt d 2^{-n_k}\le 2^{-M+d+1}r_k^{-1}$ and since $\eta(w)\ge 1$ for
$|w|\le 2^{-M+d+1}$ it follows that
\begin{equation}\label{simplelowerbd}
\eta_k*\chi_Q(y_Q)=\int r_k^d \eta(r_k(y_Q-z))\chi_Q(z) dz \ge
r_k^d2^{-n_kd}\ge 2^{-Md}.
\end{equation}

Now assume $Q\in \cQ(n_k)$ is contained in $V^N_k(x)$. For this $Q$ let
$\Omega(k,x,Q)$ be the event that $\theta_Q(\om)=1$,
but $\theta_{Q'}(\om)=0$ for all other $Q'\in \cQ(n_k)$  contained in
$V^N_k(x)$.  The probability of this event is
$\mu(\Omega(k,x,Q))= a (1-a)^{3^d 2^{(N+1)d} }$
and since $a=2^{-Nd}$  we get the uniform lower bound
\begin{equation}\label{OmkxQprob}
\mu(\Omega(k,x,Q))\ge c_d2^{-Nd}.
\end{equation}
Moreover, if $2\le l\le N$ and
$Q\subset
  Q^l_k(x)\setminus Q^{l-1}_k(x)$, $Q\subset V^N_k(x)$ then
\begin{align*}
&\int_{\Om(k,x,Q)}
\sup_{2\le j\le N}|2^{-j\sigma} I^\om_j(k,x)
|^q  d\mu  \\
&=
\int_{\Om(k,x,Q)}
2^{-lq\sigma} \sup_{y\in  Q^l_k(x)\setminus Q^{l-1}_k(x)}\big|
\eta_k*\chi_Q(y)\big|^q  d\mu
 \\
&\ge
\mu(\Om(k,x,Q)\,
2^{-lq\sigma} |\eta_k*\chi_Q(y_Q)|^q  \,\ge \,
2^{-Nd}2^{-lq\sigma} 2^{-Mdq}.
\end{align*}
For fixed $k,x$ the events
$\Omega(k,x,Q)$ are disjoint and we can sum over $Q$.
Thus
\begin{align*}
&\int_\Om
\sup_{2\le j\le N}|2^{-j\sigma} I^\om_j(k,x)
|^q  d\mu
\ge \sum_{Q\in\cQ(n_k)
\atop Q\subset V^N_k(x) }
\int_{\Om(k,x,Q)}
\sup_{2\le j\le N}|2^{-j\sigma} I^\om_j(k,x)
|^q  d\mu
\\
&\ge c_1
\sum_{2\le l\le N}
\sum_{Q\in\cQ(n_k)
\atop
{Q\subset
  Q^l_k(x)\setminus Q^{l-1}_k(x)}  }
2^{-lq\sigma} 2^{-Nd}
\ge  c_2 \sum_{2\le l\le N} (2^{ld}-2^{(l-1)d}) 2^{-lq\sigma} 2^{-Nd}
\\&\ge c_3 2^{-Nd} \max\{2^{N(d-q\sigma)}, N\}
\end{align*}
where the constants depend only on $d$, $\sigma$  and $M$.
This proves \eqref{I-lowerbd} and
\eqref{claimoflowerbound} follows.\qed

\medskip

\section{Deterministic examples}\label{deterministic}
We return to Theorem \ref{thmpeetre} and give a nonprobabilistic proof for
the lower bound in the case where $r_k=2^{-k}$, $k>0$.
With small modifications the argument can be made to apply in the general lacunary case,
where $\inf_k r_{k+1}/r_{k}>1$, but we leave this to the reader.

Fix $M>0$ sufficiently large
and let $\eta\in \cS\cap \cE_o(1)$ be a Schwartz function such that
$\eta(x)\ge 1$ if $|x_i|\le 2^{-M}$ for $i=1,\dots, d$.
Let $\eta_k= 2^{kd}\eta(2^k\cdot)$.

We fix $N$ large and set $L=2^{Nd}$.
For $k\ge N$, let $\cZ_{k,N}^d=\{0,1,\dots,2^{k-N}-1\}^d$ and for
 $j=(j_1,\dots, j_d)\in \cZ_{k,N}^d$ we set
$$Q_{k,j}= [j_1 2^{-k+N}, j_1 2^{-k+N}+2^{-k-M}]
\times\dots\times
[j_d 2^{-k+N}, j_d 2^{-k+N}+2^{-k-M}].
$$
Denote by $\chi_{k,j}$ be the characteristic function of $Q_{k,j}$,
and let $h_k=\sum_{j\in\cZ_{k,N}^d} \chi_{k,j}$.
Let $f_k=h_k*\eta_k$ so that $f_k\in \cE_o(2^k)$.

\begin{proposition} For $0<q\le p<\infty$,
\begin{equation} \label{5point1}
\Big\|\Big(\sum_{k=N}^{2^{Nd}}|f_k|^q\Big)^{1/q}\Big\|_p \le C(p,q),
\end{equation}
and for $\sigma\le d/q$
\begin{equation}\label{deterministiclowerbound}
\Big\|\Big(\sum_{k=N}^{2^{Nd}}|\fM_{\sigma, 2^k}f_k|^q\Big)^{1/q}\Big\|_p \ge c(p,q,M)
\max\{2^{N(\frac dq-\sigma)}, N^{1/q}\}.
\end{equation}
\end{proposition}

\begin{proof} It is easy to see that
$|f_k|\le C_s (M_{HL}[|h_k|^s])^{1/s}$, for $s>0$;
see  the argument for \eqref{msestimate} in the proof of
Lemma~\ref{bernoulliupperestimate}.  Thus it suffices to prove
\eqref{5point1} with $f_k$ replaced by $h_k$. For the proof we may assume that
$p\ge q$, and in fact $p=nq$ for some integer $n$ (the intermediate
cases follow by interpolation).
Thus we have to show that the $L^1([0,1]^d)$ norm of
$(\sum_{k=N}^{2^{Nd}}h_k)^n$ has an upper bound depending only on $n$.
Since each $h_k$ is nonnegative,
this follows from
\begin{equation}\label{products}
\sum_{k_1,\dots, k_n\in [N,2^{Nd}]
\atop{ k_1\le k_2\le\dots \le  k_n}} \int\prod_{i=1}^n h_{k_i}(x)
dx \le C(n).
\end{equation}
In comparison with the random
case, we have lost independence; the correlation between $h_{k_i}$  and
$h_{k_{i+1}}$  is strongest when $k_{i+1}-k_{i}$ is small.
To estimate \eqref{products} observe that
the support of $h_k$ has measure $2^{-(N+M)d}$ and that
$$\meas \big( \bigcup_{ Q_{k_{i+1}, \nu}
\subset Q_{k_i},j}
Q_{k_{i+1},\nu}\big)= \begin{cases}
|Q_{k_i,j}| 2^{-(N+M)d} &\text{ if $k_{i+1}\ge k_i+N+M$},
\\
|Q_{k_i,j}| 2^{(k_i-k_{i+1})d}
&\text{ if $k_i\le  k_{i+1}\le k_i+N+M$}.
\end{cases}
$$
Thus
$$
\int \prod_{i=1}^n h_{k_i} (x)dx =\meas \big(\supp(\prod_{i=1}^n h_{k_i} )
\big)= 2^{-(N+M)d}\prod_{i=1}^{n-1} \max\{ 2^{(k_i-k_{i+1})d} , 2^{-(N+M)d}\}.
$$
We sum in $k_n, k_{n-1},\dots, k_1$ (each ranging over the integers in
$[N, 2^{Nd}]$) and \eqref{products} follows.

We now show  for $\sigma\le d/q$ the lower bound
\eqref{deterministiclowerbound}.
Let $y_{k,j}$ be the center of $Q_{k,j}$ and observe
that $\eta_k*h_{k,j}(y_{k,j})\ge 2^{-Md}$. Thus also
\begin{equation}\label{maintermj}
\fM_{\sigma, 2^k}[\eta_k*h_{k,j}](x)\ge  2^{-Md} (1+2^k|x-y_{k,j}|)^{-\sigma}
\text{ if }|x-y_{k,j}|\le 2^{-k+N-2}.
\end{equation}
 Moreover
\begin{equation}\label{errortermsjprime}
\big|\eta_k*\sum_{j'\neq j} h_{k,j'}(x)\big|
\le C_s \big(M_{HL}[|h_k|^s]\big)^{1/s}
\text{ if } |x-y_{k,j}|\le 2^{-k+N-2}
\end{equation}
so that the terms in \eqref{errortermsjprime} are negligible error terms.
By H\"older's inequality and \eqref{maintermj}, \eqref{errortermsjprime}
\begin{align}
&\Big(\int_{[0,1]^d} \Big(\sum_{k=N}^{2^{Nd}}|\fM_{\sigma,2^k} f_k(x)|^q\Big)^{p/q} dx\Big)^{1/p}\ge
\Big(\sum_{k=N}^{2^{Nd}}\int_{[0,1]^d}
|\fM_{\sigma,2^k} f_k(x)|^q dx\Big)^{1/q}
\notag
\\
&\ge c_M
\Big(\sum_{k=N}^{2^{Nd}}\sum_{j\in \cZ_{k,N}^d}\int\limits_{
x\in [0,1]^d\atop {|x-y_{k,j}|\le 2^{-k+N-2}} }
(1+2^k|x-y_{k,j}|)^{-\sigma q}\Big)^{1/q}
\label{firstmainexpression}
\\
& \qquad -C_M  \Big(\int \big(M_{HL}[|h_k|^s]\big)^{q/s}\Big)^{1/q}.
\label{subtractederror}
\end{align}
 The subtracted term in \eqref{subtractederror}  is uniformly bounded if $s<q$
 and one easily verifies that the
main term in \eqref{firstmainexpression} is bounded below by
$c 2^{N(d-\sigma q)/q}$ if $\sigma<d/q$ and by $N^{1/q}$ if $\sigma=d/q$.
Thus \eqref{deterministiclowerbound} follows.
\end{proof}

\medskip

\section{Proof of  Theorem \ref{thmsobolev}}\label{pfthmsobolev}

We shall first  use arguments from singular integral theory to establish the upper bounds.
Then we show the lower bounds by somewhat more technical
variants of the ideas used above to prove Theorem~\ref{thmpeetre}.

\subsection{Upper bounds}
In this section we set $L_k f=\eta_k*f$ where
 $\eta_k=2^{kd}\eta(2^k\cdot)$ and $\eta$ is a
Schwartz function whose Fourier transform is supported in
 $\{\xi:1/2\le |\xi|\le 2\}$.

It suffices to set $r=2^{2^{Nd}}$ and the claimed upper bound  follows
easily
from
\begin{equation}
\Big\|
\Big(\int_0^1 \sup_{|h|\le t} \Big|\Delta_h^m
\sum_{k=1}^{2^{Nd}}
 L_k f(x)|^q t^{-1-\si q} dt\Big)^{1/q}
\Big\|_p \lc A_N(p,q,\sigma)
\Big\|\Big(\sum_{k=1}^{2^{Nd}}2^{\sigma q k}|L_k f|^q\Big)^{1/q}\Big\|_p
\end{equation}
for $q\le p$, where
$$A_N(p,q,\sigma)= 2^{N(-\sigma+d/q)}\quad \text {if $d/p<\sigma<d/q$},
\quad A_N(p,q,d/q)= \max\{N, N^{1/q}\}.$$

The contributions for the terms with $|2^kh|\le 1$ can be dealt with by
standard arguments using Peetre's maximal function. One obtains
\begin{equation}
\Big\|
\Big(\int_0^1
\sup_{|h|\le t} \Big|
\sum_{1\le k\le 2^{Nd}\atop
2^k|h|\le 1 }
\Delta_h^m
 L_k f(x)|^q t^{-1-\si q} dt\Big)^{1/q}
\Big\|_p \lc
\Big\|\Big(\sum_{k=1}^{2^{Nd}}2^{k\sigma q }
|L_k f|^q\Big)^{1/q}\Big\|_p;
\end{equation}
it is only here that the more detailed structure of the difference operator $\Delta_h^m$
and in particular the condition $m>\sigma$ is used.
Therefore  (after a change of variable and application of the triangle inequality) matters are reduced to the inequality
\begin{equation} \label{maintermsforupperbound}
\Big\|
\Big(\int_0^1\sup_{|h|\le t} \Big|
\sum_{1\le k\le 2^{Nd}\atop
2^k|h|\ge 1 }
 L_k f(x+h)|^qt^{-1-\si q} dt\Big)^{1/q}
\Big\|_p \lc A_N(p,q,\sigma)
\Big\|\Big(\sum_{k=1}^{2^{Nd}}2^{\sigma q k}|L_k f|^q\Big)^{1/q}\Big\|_p.
\end{equation}

For $n\ge 0$ define
\begin{equation}
\label{Mndefinition}
\cM^n_{k} f_k(x)
= \sup_{|h|\le 2^{n-k+1}}|f_k(x+h)|.
\end{equation}

\begin{propositionsub}\label{Mnprop}
 Let $0<q\le p<\infty$.
Then if $f_k \in \cE(2^{k})$
\begin{equation}\label{Mnpropeq}
\big\|\{\cM^n_{k} f_k\}\big\|_{L^p(\ell^q)}
\lc 2^{nd/q}
\big\|\{f_k\}\big\|_{L^p(\ell^q)}.
\end{equation}
\end{propositionsub}

\noi{\it Remarks.}

(i) Note that $\cM^n_{k} f_k\lc 2^{nd/\rho}\fM_{\rho, 2^k}f_k$ so that
the non-endpoint  $L^p(\ell^q)$  bound with constant
$C_\eps 2^{n(d/q +\eps)} $ follows from
 Peetre's maximal theorem.

(ii) There is also an endpoint inequality when $p<q$, namely
$$\big\|\{\cM^n_{k} f_k\}\big\|_{L^p(\ell^q)}
\lc 2^{nd/p}
\big\|\{f_k\}\big\|_{L^p(\ell^q)}, \qquad 0<p\le  q<\infty.
$$
This bound is not needed here and
 can be proved using arguments in \S3 of \cite{SeStud}.

\medskip

\noi{\it Proof that  Proposition \ref{Mnprop} implies
\eqref{maintermsforupperbound}.}

Assuming $q\ge 1$ we estimate

\begin{align*}
&\Big(\int_0^1\sup_{|h|\le t} \Big|
\sum_{1\le k\le 2^{Nd}\atop
{2^k|h|\ge 1 }
}
 L_k f(x+h)\Big|^q t^{-1-\si q} dt\Big)^{1/q}
\\
&\lc\sum_{m\ge 0}\sum_{n\ge 0} \Big(\sum_{l=1}^{2^{Nd}}
2^{l\sigma q}
\sup_{2^{-l-m}\le |h|\le 2^{-l-m+2}}|L_{l+m+n}f(x+h)|^q\Big)^{1/q}
\end{align*}
and using Proposition \ref{Mnprop} we obtain
\begin{align*}
&\Big\|\Big(\sum_{l=1}^{2^{Nd}}
2^{l\sigma q}
\sup_{2^{-l-m}\le |h|\le 2^{-l-m+2}}|L_{l+m+n}f(x+h)|^q\Big)^{1/q}\Big\|_p
\\&\lc 2^{-m\sigma} 2^{n(\frac dq-\sigma)}
\Big\|\Big(\sum_{l}
2^{(l+m+n)\sigma q}
|L_{l+m+n}f|^q\Big)^{1/q}\Big\|_p.
\end{align*}

The contributions of very large parameters $n$ are negligible, but
an alternative bound is needed to quantify this. One such bound
can derived by invoking H\"older's inequality to get
\begin{align*}
&\Big\|\Big(\sum_{l=1}^{2^{Nd}}
2^{l\sigma q}
\sup_{2^{-l-m}\le |h|\le 2^{-l-m+2}}|L_{l+m+n}f(x+h)|^q\Big)^{1/q}\Big\|_p
\\
&\le 2^{Nd (1/q-1/p)}
\Big\|\Big(\sum_{l=1}^{2^{Nd}}
2^{l\sigma p}
\sup_{2^{-l-m}\le |h|\le 2^{-l-m+2}}|L_{l+m+n}f(x+h)|^p\Big)^{1/p}\Big\|_p
\\&
\lc 2^{-m\sigma}2^{n(\frac dp-\sigma)} 2^{Nd (1/q-1/p)}
\Big\|\Big(\sum_{l}
2^{(l+m+n)\sigma q}
|L_{l+m+n}f|^q\Big)^{1/q}\Big\|_p.
\end{align*}

Consequently after summing in $m,n$ we obtain
\begin{equation} \cA_{p,q,\sigma}(2^{2^{Nd}})
\lc \sum_{n\ge 0}\min\{
2^{n(\frac dp-\sigma)} 2^{Nd (1/q-1/p)}, \, 2^{n(\frac dq-\sigma)} \}.
\end{equation}
In Theorem~\ref{thmsobolev} we have the hypotheses $q<p$ and $\sigma>d/p$.
Thus the series is $O(2^{N(d/q -\sigma)})$ if $\sigma < d/q$ and
is $O(N)$ when $\sigma=d/q$.

If  $q<1$ we have to bound
the $L^p$ norm of
$$\Big(
\sum_{m\ge 0}\sum_{n\ge 0} \sum_{l=1}^{2^{Nd}}
2^{l\sigma q}
\sup_{2^{-l-m}\le |h|\le 2^{-l-m+2}}|L_{l+m+n}f(x+h)|^q
\Big)^{1/q}
$$
and now
$\cA_{p,q,\sigma}^q
(2^{2^{Nd}})\lc \sum_{n=0}^\infty\min\{
2^{nq(\frac dp-\sigma)} 2^{Nd (1-q/p)}, \, 2^{nq(\frac dq-\sigma)} \}$
 which is $O(2^{N(d -q\sigma)})$ if $\sigma < d/q$ and
$O(N)$ when $\sigma=d/q$. Thus we have shown that the upper bound
in Theorem \ref{thmsobolev} is implied by Proposition \ref {Mnprop}.

\medskip
\noi {\it Proof of Proposition \ref {Mnprop}.}

We first observe that the known arguments in Peetre's maximal inequality
yield  the assertion  for $p=q$. Indeed a small modification of the proof in
\cite{Tr1}, p. 20, shows that for $g\in \cE(1)$, $0<r<\infty$, $\rho>d/r$

\begin{equation}\label{pointwisemodpeetre1}
\sup_z\frac{|\nabla g(x+z)|}{(1+2^{-n}|z|)^\rho}\lc
\sup_z\frac{|g(x+z)|}{(1+2^{-n}|z|)^\rho}
\end{equation}
and that this can be used to obtain
\begin{equation}\label{pointwisemodpeetre2}
\sup_z\frac{|g(x+z)|}{(1+2^{-n}|z|)^\rho}
\lc 2^{nd/r}
\sum_{m=0}^\infty2^{-m(\rho-d/r)}\Big(2^{-d(n+m)}\int_{|z|\le 2^{n+m}}|g(x+z)|^r dz\Big)^{1/r}.
\end{equation}
\eqref{pointwisemodpeetre2} implies for $0<r<\infty$ the inequality
\begin{equation}\label{scalarMN}
\big\|\cM^n_{k} g\big\|_{r}
\lc 2^{nd/r}
\|g\|_{r}, \quad g\in \cE(2^{k}),
\end{equation}
first for $k=0$ and then by scaling also for general $k$.
Thus we obtain
\eqref{Mnpropeq} for $p=q$.

We now consider the  assertion for $p>q$. First observe that
the case $1\le q<p$  can be  proved by interpolation
with the $L^p(\ell^p)$ bound once the cases
 $q\le p\le 1$ and $q\le 1, p>1$ are settled. We consider these cases
in what follows and use rather standard arguments from
singular integral theory, namely the
Fefferman-Stein $\#$-function estimate (\cite{FS-Hp}) which is valid for
Banach-space valued
functions; it is applied here to
 $L^{p/q}$ functions which take values in the
Banach space $\ell^1(L^\infty)$.

In what follows the slashed integral $\displaystyle\intslash_Q$ will denote an average
over the cube $Q$. By the Fefferman-Stein theorem it suffices to bound
\begin{equation*}
\Big(
\int \Big[\sup_{Q\backepsilon \,x}
\intslash_Q \sum_k \sup_{|h|\le 2^{n-k}}\Big|
|f_k(w+h)|^q- \qsl |f_k(z+h)|^q dz\Big| dw\Big]^{p/q}  dx\Big)^{1/p}
\end{equation*}
by $2^{nd/q} \big\|\{f_k\}\big\|_{L^p(\ell^q)}$.
Since  $|a|^q-|b|^q\le |a-b|^q$ for $q\le 1$ this bound
follows from
\begin{equation} \label{goalpq}
\Big(
\int \Big[\sup_{Q\backepsilon \,x}
\intslash_Q \sum_k \sup_{|h|\le 2^{n-k}}\qsl |f_k(w+h)-f_k(z+h)|^q dz dw\Big]^{p/q}  dx\Big)^{1/p}
\lc 2^{nd/q} \big\|\{f_k\}\big\|_{L^p(\ell^q)}
\end{equation}

In what follows we denote by $\ell(Q)$ the integer $\ell$ for which the sidelength of $Q$ is in $[2^{\ell-1}, 2^{\ell})$. Moreover we let $\cR_0(Q)$
be the region of all points $x$ which have distance
$\le d 2^{\ell(Q)}$ from $Q$ and for $m>0$ let
$\cR_m(Q)$
of all points $x$  for which $d2^{\ell(Q)+m-1}
\le\dist(x,Q)<d2^{\ell(Q)+m}$. Let $\eta$ be a Schwartz function in
$\cE(2)$ whose Fourier transform is equal to $1$ in
$\{\xi:|\xi|\le 1\}$. Let $\eta_k=2^{kd}\eta(2^k\cdot)$ and $P_k f=\eta_k*f$ and observe that $P_kf_k=f_k$ if $f\in \cE(2^k)$.
The estimate
\eqref{goalpq} is a consequence  of the  following three inequalities:

\begin{equation}
\label{goalpq1}
\Big(
\int \Big[\sup_{Q\backepsilon \,x}
\intslash_Q \sum_{k>n-\ell(Q)}
 \sup_{|h|\le 2^{n-k+1}}|f_k(w+h)|^q  dw\Big]^{p/q}  dx\Big)^{1/p}
\lc 2^{nd/q} \big\|\{f_k\}\big\|_{L^p(\ell^q)},
\end{equation}
\begin{multline}
\label{goalpq2}
\Big(
\int \Big[\sup_{Q\backepsilon \,x}
\intslash_Q \sum_{k=-C_0n-\ell(Q)}^{ n-\ell(Q)}
 \sup_{|h|\le 2^{n-k+1}} |f_k(w+h)|^q  dw\Big]^{p/q}  dx\Big)^{1/p}
\\\lc (n+1)^{1/q-1/p} 2^{nd/p} \big\|\{f_k\}\big\|_{L^p(\ell^q)},
\end{multline}
and, if $C_0>d/p$,
\begin{multline}
\label{goalpq3}
\Big(\int \Big[\sup_{Q\backepsilon \,x}
\intslash_Q \sum_{k<-C_0 n-\ell(Q)}
 \sup_{|h|\le 2^{n-k}}\qsl |f_k(w+h)-f_k(z+h)|^q dz dw\Big]^{p/q}  dx\Big)^{1/p}
\\
\lc \big\|\{f_k\}\big\|_{L^p(\ell^q)}.
\end{multline}

First the estimation of the main term
\eqref{goalpq1} is rather analogous to the standard
``good-function estimate'' in Calder\'on-Zygmund theory.
We split $$f_k=P_kf_k= \sum_{m=0}^\infty
 P_k[\chi_{\cR_m(Q)}f_k]$$ and estimate for fixed $x$ and $Q$
using
\eqref{scalarMN}
\begin{multline}\label{goodfctestimateforQ}
\intslash_Q \sum_{k>n-\ell(Q)}
 \sup_{|h|\le 2^{n-k+1}}|P_k
\chi_{\cR_0(Q)}f_k](w+h)|^q
  dw \\
\lc 2^{nd}
|Q|^{-1} \sum_{k>n-\ell(Q)}\int_{\bbR^d} |P_k
\chi_{\cR_0(Q)}f_k](w)|^q  dw.
\end{multline}
It is straightforward to estimate for any $\rho$, $\rho'$
\begin{equation*}
|P_k[\chi_{\cR_0(Q)} f_k](w)|\le
C_\rho\fM_{\rho, 2^k} f_k(w) \quad\text{ if } w\in \cR_0(Q), k+\ell(Q)\ge 0
\end{equation*} and
\begin{multline*}| P_k[\chi_{\cR_0(Q)} f_k](w)|\le
C_{\rho_1,\rho_2}2^{-(\ell(Q)+k+m)\rho'} \fM_{\rho, 2^k} f_k(w)
\quad\\
\text{ if } w\in \cR_m(Q),  \,m\ge 1, k+\ell(Q)+m\ge 0.
\end{multline*}
Therefore from \eqref{goodfctestimateforQ}
\begin{equation}
\sup_{Q\backepsilon \,x}
\intslash_Q
 \sup_{|h|\le 2^{n-k+1}}|P_k[
\chi_{\cR_0(Q)}f_k](w+h)|^q
  dw
\le C 2^{nd} M_{HL}( |\fM_{\rho, 2^k} f_k|^q)(x)
\end{equation}
and we may use the $L^{p/q}$ boundedness of $M_{HL}$ and the Peetre
maximal theorem  to deduce
\begin{equation}\label{goalpq1main}
\Big\|
\Big[\sup_{Q\backepsilon \,x}
\intslash_Q \sum_{k>n-\ell(Q)}
 \sup_{|h|\le 2^{n-k+1}}|P_k
[\chi_{\cR_0(Q)}f_k](w+h)|^q
  dw\Big]^{1/q}\Big\|_p \lc 2^{nd/q}\big\|\{f_k\}\big\|_{L^p(\ell^q)}
\end{equation}
We also obtain for $m\ge 1$
\begin{equation*}
|P_k[\chi_{\cR_m(Q)} f_k](x')|\le
C_\rho 2^{-(k+\ell(Q)+m)\rho}
\fM_{\rho, 2^k} f_k(x)
\end{equation*}
if $ x\in \cR_0(Q)$, $x'\in \cR_0(Q)$,
$k+\ell(Q)+m\ge 0$.
This can be applied  to bound the expression
$ \sup_{|h|\le 2^{n-k+1}}|P_k[
\chi_{\cR_m(Q)}f_k](w+h)|$ when $w\in Q$ and $k\ge n-\ell(Q)$. We
obtain
\begin{equation} \label{goalpq1error}
\Big\|
\Big[\sup_{Q\backepsilon \,x}
\intslash_Q \sum_{k>n-\ell(Q)}
 \sup_{|h|\le 2^{n-k+1}}|P_k
[\chi_{\cR_m(Q)}f_k](w+h)|^q
  dw\Big]^{1/q}\Big\|_p \lc 2^{-m\rho}\big\|\{f_k\}\big\|_{L^p(\ell^q)}
\end{equation}
and \eqref{goalpq1} follows from \eqref{goalpq1main} and
\eqref{goalpq1error}.

By H\"older's inequality the 
`left hand side of \eqref{goalpq2} is controlled by
\begin{align}
&C\Big(
\int \Big[\sup_{Q\backepsilon \,x}
\intslash_Q
\Big(\sum_{k=-C_0 n-\ell(Q)}^{ n-\ell(Q)}
|\cM^n_k f_k(w)|^q(w)\Big)^{q/p} (1+n)^{1-q/p} dw \Big]^{p/q} dx
\Big)^{1/p}
\notag
\\
&\lc (1+n)^{1/q-1/p}
\Big(\int \Big[ M_{HL}\big((\sum_k|\cM^n_k f_k|^p)^{q/p}\big)\Big]^{p/q} dx
\Big)^{1/p}
\notag
\\
&\lc (1+n)^{1/q-1/p}
\Big(\int \sum_k|\cM^n_k f_k|^p dx
\Big)^{1/p}.
\label{MnkLplp}
\end{align}
By \eqref{scalarMN} for $r=p$ we bound
\eqref {MnkLplp} by a constant times
\begin{equation*}
 (1+n)^{1/q-1/p}2^{nd/p}
\Big(\sum_k \int | f_k|^p dx
\Big)^{1/p} \lc (1+n)^{1/q-1/p}2^{nd/p}\big\|\{f_k\}\big\|_{L^p(\ell^q)}
\end{equation*}
and \eqref{goalpq2} is proved.

Finally to see  \eqref{goalpq3}
we simply observe that
$\cM^n_k(\nabla f_k)\lc 2^k
\cM^n_k(f_k)$ and thus for $x\in Q$
$$
\intslash_Q
 \sup_{|h|\le 2^{n-k}}\qsl |f_k(w+h)-f_k(w+h+z)|^q dz dw
\lc 2^{kq} (\diam(Q))^q
|\cM^n_k f_k(x)|^q.$$
Therefore the left hand side of
\eqref{goalpq3} is bounded by
\begin{equation*}
\Big(\int \big( M_{HL}[2^{-C_0 n q}|\cM_k^n f_k|^q]\big)^{p/q}\Big)^{1/p}
\lc 2^{n(\frac dp-C_0 )}\big\|\{f_k\}\big\|_{L^p(\ell^q)}. \qed
\end{equation*}

\noi{\it Remarks.} (i)  For the sequence of dyadic radii $r_k=2^k$ consider
the maximal operators  $\fM_{\sigma, r_k}$. Proposition \ref{Mnprop} can be used to show a converse to the lower bound
in Theorem \ref{thmpeetre} in this case; {\it i.e.} if
$$\cB_{p,q,\sigma}(L)
=\sup\big\{ \big\|\big(\sum_{k=1}^L |\fM_{\sigma, 2^k}f_k|^q)^{1/q}
\big\|: \, \|\{f_k\}\|_{L^p(\ell^q)}\le 1, f_k\in \cE(2^k), k=1,\dots, L\big\}
$$
then
$\cB_{p,q,\sigma}(L)\approx L^{-\sigma+d/q}$
 if $d/p<\sigma<d/q$ and
$\cB_{p,q,d/q}(L)\approx \log^{1/q}\!L$ if $p>q$.

(ii) Proposition  \ref{Mnprop} and the preceding remark remain valid for
a general lacunary sequence ($r_{k+1}/r_k\ge \gamma>1$).
\medskip

\subsection{Lower bounds}
We shall work with the Schwartz function $\eta$ defined as in the proof of
Theorem \ref{thmpeetre} (i.e. with $\widehat \eta$ vanishing
identically outside of $\{1/2< |\xi|<1\}$
and with $|\eta(x)| \ge 1$ for $|x|\le 2^{-M+d+2}$, for suitable $M$).
We shall need a
$C^\infty$ function  $\phi$ supported in $[-2^{-2M-4},2^{-2M-4}]^d$
so that
\begin{equation}\label{lowboundpsieta}
|\phi*\eta(z)|\ge c_0(M)>0
\text{ if } |z|\le 2^{-M+d+1}
\end{equation}

Let  $R$ be a large positive integer, to be chosen later. We may assume that $R\ge 10 d(1+1/p+1/q)$.
It clearly suffices to prove the lower bound in \eqref{equivalence-r}
for $r$ of the form
\begin{equation} \label{restrictiononr}
r=2^{R2^{Nd}}
\end{equation} uniformly for all large positive integers $N$.

We let
\begin{align}\label{sparseness} &n_k= kR,\quad k=1,2,\dots,
\\
& r_k=2^{n_k-M}.
\end{align}
Set $\phi _k(x) =r_k^d \phi (r_k x)$ and $\eta_k(x) =r_k^d\eta(r_k x)$.
Specify
\begin{equation}a=2^{-Nd}=L^{-1},\end{equation}
let $g_k^\om\equiv g_k^{\om,a}$ as in \eqref{definitionofgkom}, and set
\begin{equation} \label{Gdef}
G_k^\om(x)=2^{-n_k\sigma} g_k^\om(x), \qquad G^\omega(x)=
\sum_{k=1}^{2^{Nd}}
G_k^\om(x).
\end{equation}

We need the following
estimates for convolutions with the functions $g_k^\om$ and $G^\om$.


\begin{lemmasub}\label{lemmaupperbounds1}
(i)
 Let $H$ be a Schwartz function so that $\int H(x) x^\alpha dx=0$ for all
multiindices $\alpha$ with $\max_i|\alpha|_i\le N_0$.
Let $H_\ell=2^{\ell d}H(2^\ell\cdot)$. Then
\begin{equation}\label{convolutionwithHell}
|H_\ell*g_l^\om(x)|\lc 2^{-|\ell-n_l|(N_0-d/s)}
\big(M( (\sum_{Q\in\cQ(n_l)} \theta_Q(\om)\chi_Q)^s)(x)\big)^{1/s}.
\end{equation}

(ii)
 For $0<p,q <\infty$
\begin{equation}
\Big(\int_\Om\|G^\omega\|^p_{F^{p}_{\si q}}d\mu(\omega)\Big)^{1/p}\le C_1
\end{equation}
and
\begin{equation}
\Big(\int_\Omega
\int_{\bbR^d}\Big[\int_0^1\sup_{|h|\le t}
\sum_{l: t2^{n_l} \le 1}
\big|
\Delta_h^m G_l^\om(x)
\big|^q t^{-1-\si q} dt\Big]^{p/q} dxd\mu(\omega) \Big)^{1/p}
\le C_2.
\end{equation}
Moreover
\begin{equation}
\Big(\int_\Omega
\int_{\bbR^d}\Big[\int_0^1\Big(
\sum_{l: t2^{n_l} \ge 1}
|G_l^\om(x)|\Big)^q t^{-1-\si q} dt\Big]^{p/q} dxd\mu(\omega) \Big)^{1/p}
\le C_3.
\end{equation}
Here $C_1$, $C_2$, $C_3$ depend only on  $p$,$q$, $m$ and $d$.
\end{lemmasub}

\begin{proof}
\eqref{convolutionwithHell} is straightforward, \cf. the reasoning for    inequality
\eqref{msestimate}.
The other assertions follow in a  straightforward manner from the
basic estimates \eqref{msestimate} and \eqref{convolutionwithHell},
a suitable application of Minkowski's inequality, and
Lemma   \ref{bernoulliupperestimate}.
\end{proof}

Let $\tDel^m_h $ be defined by
\begin{equation}\label{definetildedelta}
\tDel^m_h f(x)=\Delta^m_h f(x)-(-1)^m f(x)=\sum_{\nu=1}^m (-1)^{m-\nu}
\binom{m}{\nu}
f(x+\nu h)
\end{equation}
and let
$$I_{k,j}=[2^{-n_k+j+2d}, 2^{-n_k+j+2d+1}].$$
In view of Lemma \ref{lemmaupperbounds1} and H\"older's inequality   on
$\Omega \times [0,1]^d$ (with $p/q\ge 1$),
in order to prove the lower bounds in
\eqref{equivalence-r}, \eqref{equivalence-r-2}, \eqref{inequality-r}
it suffices to prove that
\begin{multline}\label{lb1}
\Big(\int_\Omega
\int_{[\tfrac 14,\tfrac 34]^d}
\sum_k
\sum_{j=d+M}^{N-M-d}\int_{I_{k,j}}
\sup_{|h|\le mt}\Big|
\sum_{l\ge k}
\tDel_{\tfrac hm}^m G_l^\om(x)
\Big|^q t^{-1-\si q} dt
 dxd\mu
 \Big)^{1/q}
\\
\ge c_0 \max\{2^{N(\frac dq-\si)}, N^{1/q}\}
\end{multline}
for some $c_0>0$.

Let
\begin{equation}\label{defineGammaQ}
\Gamma^{l,Q}_{k,m}(x,h)
:=\int\phi _{k}(y)
\tDel_{\tfrac {h-y}m}^m \eta_l*\chi_Q(x)dy, \qquad Q\in \cQ(n_l)
\end{equation}
and
\begin{equation}\label{defineGamma}
\Gamma^l_{k,m}(x,h,\om) :=
\int\phi _{k}(y)
\tDel_{\tfrac {h-y}m}^m G_l^\om(x)dy = 2^{-n_l\sigma}\sum_{Q\in \cQ(n_l)}
\theta_Q(\om)
\Gamma^{l,Q}_{k,m}(x,h).
\end{equation}
We use the elementary inequality
$$|\phi _{k}*a(h)|\le C \sup_{|h-u|\le 2^{-n_{k}-2}}|a(u)|
$$
to deduce that \eqref{lb1} follows from the existence of a constant $c_1>0$ such that
\begin{multline} \label{lb2}
\Big(\int_\Omega
\int_{[\tfrac 14,\tfrac 34]^d}
\sum_k \sum_{j=M+d}^{N-M-d}
2^{(n_k-j)\sigma q}\sup_{|h|\le m2^{-n_k+j+d}}\Big|
\sum_{l\ge k}
\Gamma^l_{k,m}(x,h,\om)
\Big|^q dxd\mu \Big)^{1/q}
\\
\ge c_1 \max\{ 2^{N(\frac dq-\si)}, N^{1/q}\}.
\end{multline}

We show now  that the only relevant terms in \eqref{lb2}
are those with $l=k$; it is here where we have to choose $R$ sufficiently large.

\begin{lemmasub}\label{lemmaupperbounds2}
For $0<q\le p <\infty$,
\begin{multline}\label{lemmaupperbounds2equation}
\Big(\int_\Omega
\int
\sum_k \sum_{j=M+d}^{N-M-d}
2^{(n_k-j)\sigma q}\sup_{|h|\le m2^{-n_k+j+d}}\Big|
\sum_{l\ge k+1}
\Gamma^l_{k,m}(x,h,\om)
\Big|^q  dxd\mu \Big)^{1/q}
\\ \le C_4 2^{-R} \max\{ 2^{N(\frac dq-\sigma)}, N^{1/q}\}
\end{multline}
\end{lemmasub}
\begin{proof}
Here we use the cancellation inherent in $\Psi_Q:=\eta_l*\chi_Q$,
stemming from the vanishing of $\widehat{\eta}(\xi)$ near $\xi=0$.
Define $y_{Q,\nu}\equiv
y_{Q,\nu}(x,h)=h+(x-x_Q)m/\nu$
where $x_Q$  is the center of $Q$.
 Now by \eqref{definetildedelta}
and Taylor's formula
\begin{align*}
&\Gamma^{l,Q}_{k,m}(x,h)=\sum_{\nu=1}^m c_\nu
\int \phi _k(y) \Psi_Q(x+\tfrac \nu m (h-y)) dy
\\
&=\sum_{\nu=1}^m c_\nu
\int \int_0^1\frac{(1-s)^{N_0-1}}{N_0!}\inn{y-y_{Q,\nu}}{\nabla}^{N_0}
\phi _k(y_Q+s(y-y_{Q,\nu})) ds\,
 \Psi_Q(x+\tfrac \nu m(h-y)) \,dy
\end{align*}
where  $c_\nu=(-1)^{m-\nu}\binom{m}{\nu}$,  $Q\in \cQ(n_l)$.
Let
\begin{equation}
\label{Foml}
 F^\omega_l(x) =\sum_{Q\in \cQ(n_l)} \theta_Q(\om)
\eta_l*\chi_Q(x);
\end{equation}
then we see that for $l>k$
\begin{equation}\label{Gammakmlmaxest}
\sup_{|h|\le m 2^{-n_k+j+d}} |\Gamma_{k,m}^l(x,h,\om)|
\lc 2^{-n_l\sigma q} 2^{(n_k-n_l)N_0} \sup_{|y|\le 2^{-n_k+j+d}}
|F^\om_l(x+y)|
\end{equation}
By \eqref{scalarMN}
\begin{equation}\label{Fomlestimate}
\Big\|\sup_{|y|\le 2^{-n_k+j+d}}
|F^\om_l(\cdot+y)|\Big\|_q \lc 2^{(n_l-n_k+j)d/q} \|F^\om_l\|_q
\end{equation}
and therefore the left hand side of
\eqref{lemmaupperbounds2equation} is dominated by
\begin{align*}
&\sum_{s=1}^\infty\Big(
\int_\Omega
\int
\sum_{k=1}^{2^{Nd}} \sum_{j=M+d}^{N-M-d}
2^{(n_k-j)\sigma q}\sup_{|h|\le m2^{-n_k+j+d}}\Big|
\Gamma^{k+s}_{k,m}(x,h,\om)
\Big|^q  dxd\mu \Big)^{1/q}
\\
&\lc \sum_{s=1}^\infty\Big(
\int_\Omega
\int
\sum_{k=1}^{2^{Nd}} \sum_{j=M+d}^{N-M-d}
2^{(n_k-n_{k+s}-j)\sigma q}
2^{(n_k-n_{k+s})N_0 q} \sup_{|y|\le 2^{-n_k+j+d}}
|F^\om_{k+s}(x+y)|^q
 dxd\mu \Big)^{1/q}
\\
&\lc \sum_{s=1}^\infty 2^{-Rs ( N_0+\sigma-d/q)}
 \Big(
\sum_{k=1}^{2^{Nd}} \sum_{j=M+d}^{N-M-d} 2^{-j(\sigma q-d)}
\int_\Omega
\|F^\om_{k+s}\|_q^q d\mu \Big)^{1/q}
\\
&\lc 2^{-R} \max \{2^{(-\sigma+d/q)N}, N^{1/q} \},
\end{align*}
by  \eqref{Gammakmlmaxest}, \eqref{Fomlestimate},
Minkowski's inequality and Corollary \ref{bernoulliupperestimateeta}.

If  $q<1$ we use the $\ell^q$ triangle inequality  in place of
Minkowski's inequality and we have to bound
$$\Big(\sum_{s=1}^\infty
\int_\Omega
\int
\sum_{k=1}^{2^{Nd}} \sum_{j=M+d}^{N-M-d}
2^{(n_k-n_{k+s}-j)\sigma q}
2^{(n_k-n_{k+s})N_0 q} \sup_{|y|\le 2^{-n_k+j+d}}
|F^\om_{k+s}(x+y)|^q
 dxd\mu \Big)^{1/q};
$$
the result is the same.
\end{proof}

Given Lemma \ref{lemmaupperbounds2} the
lower bound \eqref{lb2} follows from a corresponding lower bound
for the expression only involving the
$\Gamma^k_{k,m}(x,h,\om)$ and
it remains to show for $\sigma\le d/q$:
\begin{lemmasub}\label{sobolevlowerbound}
\begin{multline}\label{lb3}
\Big( \sum_k
\sum_{j=M+d}^{N-M-d}2^{-j\sigma q}
\int_{[\frac 14,\frac 34]^d} \int_\Omega
\sup_{|h|\le m2^{-n_k+j+d}}\Big|
2^{n_k\sigma}\Gamma^k_{k,m}(x,h,\om)
\Big|^q  d\mu dx\Big)^{1/q}
\\ \ge c_2 \max\{2^{N(\frac dq-\si)}, N^{1/q}\}
\end{multline}
for some $c_2>0$.
\end{lemmasub}

\begin{proof}
In what follows we fix $x\in [1/4,3/4]^d$ and $1\le k\le 2^{Nd}$.
As in the proof of Theorem \ref{thmpeetre}
define $V^N_k(x)$ to be the union of all dyadic cubes of sidelength
$2^{-n_k+N+1}$ whose boundaries intersect the boundary of $Q_k^{N+1}(x)$.
Let $\cV^N_k(x)$ be the set of all  $Q\in \cQ(n_k)$ that are
 contained in the closure of  $V^N_k(x)$. 
Denote by 
$\Omega(k,x,Q)$  the event  that $\theta_Q(\om)=1$,
but $\theta_{Q'}(\om)=0$ for all $Q'\in\cV^N_k(x)\setminus\{Q\}$.
For the probability of this event there is the  lower bound
$\mu(\Omega(k,x,Q))\ge c_d 2^{-Nd}$; see \eqref{OmkxQprob}.

Now let $\cW(k,j,x)$ be the set of all cubes $Q\in\cQ(n_k)$ for which
$$2^{-n_k+j}\le \dist(x, Q)\le 2^{-n_k+j+1}.$$
For $Q\in \cW(k,j,x)$ denote by $y_Q$ the center of $Q$ and set
$h_{Q,x}=y_Q-x$ so that $|h_{Q,x}|\lc 2^{-n_k+j+1}.$

Thus the left hand side of
\eqref{lb3} is bounded below by
\begin{equation}\label{decomposeintoevents}
c\Big( \sum_k
\sum_{j=M+d}^{N-M-d}2^{-j\sigma q}
\int_{[\frac 14,\frac 34]^d}
\sum_{Q\in \cW(k,j,x)} \int_{\Omega(k,x,Q)}
\Big|
2^{n_k\sigma}\Gamma^k_{k,m}(x,h_{Q,x},\om)
\Big|^q  d\mu dx\Big)^{1/q}.
\end{equation}

For  $Q\in \cW(k,j,x)$
 and $\omega\in \Omega(k,x,Q)$ we decompose further
$$
2^{n_k\sigma}\Gamma^k_{k,m}(x,h_{Q,x},\om)= \sum_{\nu=1}^m (-1)^{m-\nu}\binom{m}{\nu} I_{\nu}^\om(k,x,Q)  + II^\om(k,x)
$$
where
\begin{equation}
I_{\nu}^\om (k,x,Q)=\int \phi _k(y) \eta_k*\chi_Q(x+\frac{\nu}m (h_{Q,x}-y)) dy
\end{equation}
and
\begin{equation}
II^\om (k,x,Q)
= \sum_{Q'\in \cQ(n_k)\atop
{Q'\notin\cV^N_k(x)} }\theta_{Q'}(\om)
\int\phi _{k}(y)
\tDel_{\tfrac {h(Q,x)-y}m}^m \eta_l*\chi_{Q'}(x)dy.
\end{equation}
We prove a lower bound for $I_m^\om$ and upper bounds for $II^\om$ and $I_\nu^\om$, $\nu\le m-1$.

Notice that for $\omega\in \Omega(k,x,Q)$
\begin{align*}
I_{m}^\om (k,x,Q)&=\int \phi _k(y) \eta_k*\chi_Q(x+h_{Q,x}-y) dy
\\&=\int r_k^d (\eta*\phi )(r_k(y_Q-z))\chi_Q(z) dz
\end{align*}
and  since $|r_k(y_Q-z)|\le \sqrt d 2^{-n_k}\cdot 2^{n_k-M}$ for $z\in Q$ it follows from
\eqref{lowboundpsieta} that
\begin{equation}\label{lbImom}
|I_{m}^\om (k,x,Q)|
 \ge c(M) r_k^d 2^{-n_k d} \ge c'(M) \qquad\text {if }
\omega\in \Omega(k,x,Q), \, Q\in \cW(k,j,x).
\end{equation}
Since $\mu(\Omega(k,x,Q))\ge c 2^{-Nd}$
and $\card (\cW(k,j,x))
\approx 2^{jd}$,
\begin{align}
&\Big( \sum_{k=1}^{2^{Nd}}
\sum_{j=M+d}^{N-M-d}2^{-j\sigma q}
\int_{[\frac 14,\frac 34]^d}
\sum_{Q\in \cW(k,j,x)}
 \int_{\Omega(k,x,Q)}
|I_{m}^\om (k,x,Q)|^q
 d\mu dx\Big)^{1/q}
\notag
\\
&\gc\Big(\int_{[\frac 14,\frac 34]^d}
 \sum_{k=1}^{2^{Nd}} 2^{-Nd}
\sum_{j=M+d}^{N-M-d}2^{j(d-\sigma q)}
\Big)^{1/q} \gc c_q \big(\max \{2^{N(d-\sigma q)}, N\}\big)^{1/q}.
\label{lowerboundImom}
\end{align}
Next notice that for
$\omega\in \Omega(k,x,Q)$,
$Q\in \cW(k,j,x)$, $y\in \supp \phi _k$ and $\nu\le m-1$,
$$|x+(h_{Q,x}-y)\nu/m -y_Q|\ge |x-y_Q|(1-\nu/m)-\nu/m|y| \gc 2^{-n_k +j}$$
which shows that
$|\eta_k*\chi_Q(x+\frac{\nu}m (h_{Q,x}-y))| \le C_\rho 2^{-j\rho}$
for all $\rho$ and consequently we get the estimate
\begin{equation*}
|I_{\nu}^\om (k,x,Q)|\le
 C_{M,\rho} \fM_{\rho, 2^{-n_k}}[h_k^\om],
\quad \nu\le m-1.
\end{equation*}
Similarly for  $II^\om (k,x,Q)$ we can argue as for the corresponding term
in the proof of Theorem \ref{thmpeetre}
and see that
\begin{equation*}
\sup_{2\le j\le N}\big|II^\om (k,x)
|\le C_{M,\rho} \fM_{\rho, 2^{-n_k}}
[h_{k}^\om] \quad\text{ if }\om\in \Omega(k,x,Q)
\end{equation*}
for any $\rho>0$. Thus
\begin{align}
&\Big( \sum_{k=1}^{2^{Nd}}
\sum_{j=M+d}^{N-M-d}2^{-j\sigma q}
\int_{[\frac 14,\frac 34]^d}
\sum_{Q\in \cW(k,j,x)}
 \int_{\Omega(k,x,Q)}
\Big[\sum_{\nu=1}^{m-1}
|I_{\nu}^\om (k,x,Q)|
+|II^\om (k,x,Q)|\Big]^q
 d\mu dx\Big)^{1/q}
\notag
\\&
\lc\Big( \sum_{k=1}^{2^{Nd}}
\int_\Omega \int
\big(\fM_{\rho, 2^{-n_k}}
h_{k}^\om\big)^q dx
d\mu\Big)^{1/q} \le C
\label{InuIIerrors}
\end{align}
by Lemma \ref{hk}.
We combine the estimates
\eqref{lowerboundImom} and
\eqref{InuIIerrors} to see that the expression \eqref{lb3}
is bounded below by $\max\{ 2^{N(d-\sigma/q)}, N^{1/q}\}$.
\end{proof}

We now combine the various estimates and note that the lower
bound in Lemma \ref{sobolevlowerbound} is independent of $R$.
Thus if $R$ is chosen to be sufficiently large, the upper bounds in
\eqref{lemmaupperbounds2equation} can be absorbed by the lower bound
\eqref{lb3}, and \eqref{lb2} consequently follows.
All told, we have shown the lower bound $\cA_{p,q,\sigma}(2^{R2^{Nd}})\gc
\max\{ 2^{N(d-\sigma/q)}, N^{1/q}\}$, which implies the asserted
lower bound for large $r= 2^{R2^{Nd}}$. \qed

\medskip

\noi{\it Remark.} One can also consider the more regular variant
\begin{equation*}
\fD_{s,m}^{\si,q} f(x)= \Big(\int_0^1 \Big[ \intslash_{|h|\le
t}
|\Delta_h^m f(x)|^s  dh \Big]^{q/s}t^{-1-\si q}
dt\Big)^{1/q}.
\end{equation*}
where the slashed integral denotes the average over the ball $\{h:|h|\le t\}$.
 Then $\|f\|_\fpqs\approx
\|f\|_p+\|\fD_{s,m}^{\si,q} f\|_p$ provided that
$\si> \max\{0, d(1/p-1/s), d(1/q-1/s)\}$. A modification of our argument
shows that
the characterization fails when $\si\le d(1/q-1/s)$.

\section{Proof of Theorem \ref{oscillatorytheorem}}\label{pfoscillatory}
Let $\eta$ be a Schwartz function in $\cE(2)$ such that $\widehat
\eta$ vanishes identically in a neighborhood of the origin and
$\widehat \eta(\xi)=1 $ if $2^{-1/2}\le |\xi|\le 2^{1/2}$.
In what follows we fix $q<p\le 2$ and assume that $0<\gamma<1$ and  that
$b=\gamma d (1/p-1/2)$.
Define for $k>1$ the operator $T_k$ by
\begin{equation}
\widehat {T_k f}(\xi)=
e^{i|\xi|^\gamma}\widehat \eta(2^{-k}\xi)\widehat  f(\xi).
\end{equation}
It is easy to see, using \eqref{vectorvaluedmaxfct},
that the statement of the Theorem is equivalent with the statement that the
best constant $\cA_L$ in the inequality
\begin{equation}
\Big\|\Big(\sum_{k=1}^L| 2^{-k b} T_k f_k|^q\Big)^{1/q}\Big\|_p\le \cA_L
\Big\|\Big(\sum_{k=1}^L|f_k|^q\Big)^{1/q}\Big\|_p ,\qquad \text{with } f_k\in
\cE(2^{k+1}),
\end{equation}
satisfies $\cA_L\approx L^{1/q-1/p}$.
As the operators $2^{-k b} T_k$ map
$L^p\cap \cE(2^{k+1})$ to $L^p$, with  bounds uniform in $k>0$
({\it cf.}\ \cite{FS-Hp}),
 the upper bound
$\cA_L\lc  L^{1/q-1/p}$ is immediate by H\"older's inequality and the embedding $\ell^q\subset\ell^p$.
In what follows we prove the lower bound.

We use a variant of the  random construction of \S\ref{sectionrandom} and define $\theta_{Q,a}$ and
$h_{k}^{\om,a}$ as  in \eqref{hk};
however we  now let $a$  depend  on $k$ and require that
\begin{equation}\label{defofa}
a_k= 2^{-k\gamma d }.
\end{equation}
Also let $\widetilde \eta$ be a Schwartz function in $\cE(2)$
whose Fourier transform equals $1$ on the support of $\widehat \eta$.
Define  (using the notation in \S\ref{sectionrandom} with $n_k= k$)
\begin{equation}
f_k^\om(x)= \beta_k \sum_{Q\in \cQ(k)}
\theta_{Q, a_k}(\om)
\widetilde\eta(2^k(x-x_Q))
\end{equation}
where $x_Q$ is the center of $Q$ and
\begin{equation}\label{betakakrelation}
 \beta_k=a_k^{-1/p}.
\end{equation}
We claim
\begin{equation}\label{fkoscestimate}
\Big(\int_\Omega\Big\|\Big(\sum_{k=1}^L|f_k^\om|^q\Big)^{1/q}\Big\|_p^p
d\mu(\omega) \Big)^{1/p} \lc L^{1/p};
\end{equation}
this inequality will use only  \eqref{betakakrelation}   and the fact
that  the $\beta_k$ increase at least in a geometric  progression;
the  specific choice \eqref{defofa}
 is not yet needed.

A straightforward estimate yields
$
|f_k^\om(x)|\le C_s\big(M_{HL}(|\beta_k h_k^{\omega,a_k}|^s)\big)^{1/s}
$
for all $s>0$
and therefore it suffices to prove
\eqref{fkoscestimate} with $f_k^\om$ replaced by
$\beta_k h_k^{\omega,a_k}$.

Now since $h_k^{\om,a_k}$ takes only values $1$ and $0$ and the $\beta_k$
increase at least geometrically we see that for all $x$
$$
\Big(\sum_{k=1}^L|\beta_k h_k^{\om,a_k}(x)|^q\Big)^{1/q}
\le C_\gamma
\sup_{1\le k\le L}|\beta_k h_k^{\om,a_k}(x)|,
$$
with a finite constant $C_\gamma$ independent of $x$.
After replacing the supremum by an $\ell^p$ norm we see that the left hand
side of \eqref{fkoscestimate} can be estimated by
\begin{equation*}\Big(\int_{[0,1]^d} \sum_{k=1}^L\int_{\Omega}
\big|\beta_k h_k^{\om,a_k}\big|^p
d\mu(\omega) dx\Big)^{1/p}
\le C
\Big(\int_{[0,1]^d} \sum_{k=1}^L \beta_k^p  a_k dx\Big)^{1/p} =
C L^{1/p}.
\end{equation*}

It remains to show the lower bound
\begin{equation}\label{lowerboundoscfortheLqlqexpression}
\Big(\int_\Omega\Big\|\Big(\sum_{k=1}^L|T_kf_k^\om|^q\Big)^{1/q}\Big
\|_p^p
d\mu \Big)^{1/p} \ge c L^{1/q}.
\end{equation}
Now let $K_k$ be the convolution kernel of $T_k$. Then
$$T_k f^\omega_k(x)= \beta_k
\sum_{Q\in \cQ(k)}\theta_{Q,a_k}(\om) 2^{-kd} K_k(x-x_Q).$$
A stationary phase calculation shows that
for suitable $\eps_1>0$ there is the uniform estimate for large $k$
\begin{equation}\label{lowerboundforKk}
|K_k(x)|\ge 2^{k(d-d\gamma/2 )}\quad \text{ if }
(1-\eps_1)
2^{-k(1-\gamma)}
\le |x|
\le (1+\eps_1)
2^{-k(1-\gamma)};
\end{equation}
moreover for any $\rho<\infty$
\begin{equation}\label{upperboundforKk}
|K_k(x)|\le C_\rho2^{kd} (2^{k}|x|)^{-\rho}
\text{ if } |x|\ge  B 2^{-k(1-\gamma)}
\end{equation}
for suitable $B$ ($\ge 2$); this is seen by using integration by parts
for the oscillatory integral,  which has  a nonstationary phase
when $|x|\ge  B 2^{-k(1-\gamma)}$.
Now apply H\"older's inequality (as in all previous examples):
\begin{multline*}\Big(\int_\Om\int_{[0,1]^d} \Big(\sum_{k=1}^L
\big|2^{-kb} T_kf_k^\om\big|^q\Big)^{p/q}
dxd\mu \Big)^{1/p}
\\ \ge
\Big(\int_{[\frac 14, \frac 34]^d}\int_\Om\sum_{k=1}^L \Big|
2^{-k(d+b)} \beta_k\sum_{Q\in\cQ(k)} \theta_{Q,a_k}(\om)  K_k(x-x_Q)\Big|^q
d\mu dx \Big)^{1/q}.
\end{multline*}

Fix $x\in [1/4,3/4]^d$ and let
$\cV_{k,\gamma}(x)$ be the set
of all cubes $Q\in\cQ(k)$ whose distance to $x$ is $\le  C_2  2^{-k(1-\gamma)}$ where
$C_2\gg  B$ for a sufficiently large constant $B$.
Let $\Omega(k,x,Q)$
 be the event that
$\theta_Q(\om)=1$ but $\theta_{Q'}(\om)=0$ for all $Q'\in
\cV_{k,\gamma}(x)\setminus\{Q\}$.
The probability of this event satisfies
$$\mu(\Omega(k,x,Q))\ge a_k (1-a_k)^{\card(\cV_{k,\gamma}(x))-1}
\ge c a_k,
$$ by our choice \eqref{defofa}.


By the upper bound \eqref{upperboundforKk} we get
\begin{align*}
&\int_\Om\Big|
2^{-k(d+b)} \beta_k \sum_{Q\in \cQ(k)\setminus \cV_{k,\gamma}(x)} \theta_{Q,a_k}(\om)
K_k(x-x_Q)\Big|^q
d\mu(\omega)
 \\&\lc 2^{kd}
\int_{|z-x|\ge 2^{-k(1-\gamma)}}
(\beta_k 2^{-k (b+d)} )^q (2^k|z-x|)^{-\rho q}\,dz
\le \beta_k^q 2^{-kbq} 2^{k\gamma(d-\rho q)}
\end{align*}
provided that $\rho>d/q$.
Consequently by choosing $\rho$ large enough we find that
\begin{equation}\label{upperboundforoscerror}
\Big(\int_{[\frac 14, \frac 34]^d}\int_\Om\sum_{k=1}^L \Big|
2^{-k(d+b)} \beta_k\sum_{Q\in \cQ(k)\setminus \cV_{k,\gamma}(x)} \theta_{Q,a_k}(\om)
K_k(x-x_Q)\Big|^q
d\mu dx \Big)^{1/q} \le C
 \end{equation}
uniformly in $L$.

By \eqref{upperboundforoscerror} we may estimate
\begin{align}
&\Big(\int_{[\frac 14, \frac 34]^d}\int_\Om\sum_{k=1}^L \Big|
2^{-k(d+b)} \beta_k \sum_Q \theta_{Q,a_k}(\om)  K_k(x-x_Q)\Big|^q
d\mu dx \Big)^{1/q}
\notag
\\&\ge
\Big(\int_{[\frac 14, \frac 34]^d}\sum_{k=1}^L
\int_\Om\Big|
2^{-k(d+b)} \beta_k\sum_{Q'\in \cV_{k,\gamma}(x)} \theta_{Q',a_k}(\om)  K_k(x-x_{Q'})\Big|^q
d\mu dx \Big)^{1/q} -C
\notag
\end{align}
and the main term is
\begin{align}&\ge
\Big(\int_{[\frac 14, \frac 34]^d}\sum_{k=1}^L
\sum_{Q\in \cV_{k,\gamma}(x)}
\int_{\Om(k,x,Q)}\Big|
2^{-k(d+b)} \beta_k \sum_{Q'\in \cV_{k,\gamma}(x)} \theta_{Q',a_k}(\om)  K_k(x-x_{Q'})\Big|^q
d\mu dx \Big)^{1/q}
\notag
\\&=
\Big(\int_{[\frac 14, \frac 34]^d}\sum_{k=1}^L
\sum_{Q\in \cV_{k,\gamma}(x)}
\int_{\Om(k,x,Q)}\Big|
2^{-k(d+b)}   \beta_k K_k(x-x_{Q})\Big|^q
d\mu dx \Big)^{1/q}.
\label{lasttechnicalest}
\end{align}
Now let $\cW_{k,\gamma}(x)$ be the family of all cubes in $\cQ(k)$
which are contained in the set
$\{y: (1-\eps_1)2^{-k(1-\gamma)}\le |x-y|\le
(1+\eps_1)2^{-k(1-\gamma)}\}$.
These cubes are also  in $\cV_{k,\gamma}(x)$ and if
 $Q\in\cW_{k,\gamma}(x)$
then we may use the lower bound
\eqref{lowerboundforKk} for the term $K_k(x-x_Q)$.
Note also that
$\card( \cW_{k,\gamma}(x))\gc 2^{kd\gamma}\approx a_k^{-1} $
for large $k$.

Thus the term \eqref{lasttechnicalest} is bounded below for large $L$  by
\begin{multline*}
\Big(\int_{[\frac 14, \frac 34]^d}\sum_{k=1}^L
\sum_{Q\in \cW_{k,\gamma}(x)}
\int_{\Om(k,x,Q)}\big[
2^{-k(d+b)}  \beta_k2^{k(d-\gamma d/2)}\big]^q d\mu\, dx\Big)^{1/q}
\\
\ge c
\Big(\int_{[\frac 14, \frac 34]^d}\sum_{k=C}^L
\card(\cW_{k,\gamma}(x)) a_k dx\Big)^{1/q} \gc L^{1/q}
\end{multline*}
and consequently we obtain
\eqref{lowerboundoscfortheLqlqexpression}.\qed

\medskip

\noi{\it Remark:} The case $\gamma=1$ which is relevant
for the wave equation is an exceptional case
(see \cite{Miyachi}, \cite{Peral}),
as the critical $b$ is given by $b=(d-1)(1/p-1/2)$, $1<p\le 2$.
However if these parameters are chosen in Theorem \ref{oscillatorytheorem}
then  a modification of the above argument, with $a_k=2^{-k(d-1)}$,  
shows that
\eqref{oscequiv} remains valid.

\end{document}